\documentclass[12pt]{article}
\usepackage{amsmath,amssymb,amsthm,color,hhline,euscript}

\newcommand{\fr}[2]{\frac{#1}{#2}}
\newcommand{\dfr}[2]{\dfrac{#1}{#2}}
\newcommand{\cd}{\cdot}
\newcommand{\cds}{\cdots}

\renewcommand{\l}{\left}
\renewcommand{\r}{\right}

\newcommand{\vsv}{\vspace{5mm}}
\newcommand{\vsb}{\vspace{2mm}}

\newcommand{\q}{\quad}
\newcommand{\qq}{\qquad}

\newcommand{\maru}[1]{{\ooalign{\hfil#1\/\hfil\crcr
\raise.167ex\hbox{\mathhexbox20D}}}}

\newcommand{\ruby}[2]{%
 \leavevmode
 \setbox0=\hbox{#1}%
 \setbox1=\hbox{\tiny #2}%
 \ifdim\wd0>\wd1 \dimen0=\wd0 \else \dimen0=\wd1 \fi
 \hbox{%
   \kanjiskip=0pt plus 2fil
   \xkanjiskip=0pt plus 2fil
   \vbox{%
     \hbox to \dimen0{%
       \tiny \hfil#2\hfil}%
     \nointerlineskip
     \hbox to \dimen0{\mathstrut\hfil#1\hfil}}}}

\newcommand{\la}{\langle}
\newcommand{\ra}{\rangle}

\newcommand{\abs}[1]{\lvert{#1}\rvert}

\DeclareMathOperator*{\tensor}{\otimes}

\newcommand{\Z}{\mathbb{Z}}
\newcommand{\C}{\mathbb{C}}
\newcommand{\R}{\mathbb{R}}
\newcommand{\N}{\mathbb{N}}
\newcommand{\Q}{\mathbb{Q}}

\newcommand{\vir}{\mathrm{Vir}}
\newcommand{\aut}{\mathrm{Aut}}

\renewcommand{\hom}{\mathrm{Hom}}

\newcommand{\id}{\mathrm{id}}
\newcommand{\ind}{\mathrm{Ind}}

\newcommand{\hf}{\fr{1}{2}}
\newcommand{\Span}{\mathrm{Span}}

\newcommand{\w}{\omega}
\newcommand{\vacuum}{\mathrm{1\hspace{-3.2pt}l}}
\newcommand{\vac}{\vacuum}

\newcommand{\supp}{\mathrm{supp}}

\newcommand{\sfr}[2]{\leavevmode\kern-.1em
  \raise.5ex\hbox{\the\scriptfont0 #1}\kern-.1em
  /\kern-.15em\lower.25ex\hbox{\the\scriptfont0 #2}}

\newcommand{\shf}{\sfr{1}{2}}

\makeatletter
\@addtoreset{equation}{section}

\theoremstyle{plain}
\newtheorem{thm}{Theorem}[section]
\newtheorem{prop}[thm]{Proposition}
\newtheorem{lem}[thm]{Lemma}
\newtheorem{cor}[thm]{Corollary}

\theoremstyle{definition}
\newtheorem{df}[thm]{Definition}

\theoremstyle{remark}
\newtheorem{rem}[thm]{Remark}

\newcommand{\pf}{\noindent {\bf Proof:}\q }

\newcommand{\com}{\mathrm{Com}}

\newcommand{\om}{\omega}
\newcommand{\be}{\beta}
\newcommand{\al}{\alpha}
\newcommand{\eL}{\EuScript{L}}

\title{Ising vectors and automorphism groups of commutant subalgebras
related to root systems}
\date{}
\author{Ching Hung Lam\footnote{Supported by NSC grant 93-2115-M-006-012 of Taiwan
and National Center for Theoretical Sciences, Taiwan}
 \vsb\\
{\small \sl Department of Mathematics, National Cheng Kung University,}
\\
{\small \sl Tainan, Taiwan 701}
\\
{\small ${}^*$e-mail: {\tt chlam@mail.ncku.edu.tw}}
\vsb\\
Shinya Sakuma\footnote{Supported by JSPS Research Fellowships for Young
Scientists.}\, and
Hiroshi Yamauchi\footnote{Supported by JSPS Research Fellowships for Young Scientists.}
\vsb\\
{\small \sl Graduate School of Mathematical Sciences, The University of Tokyo,}\\
{\small \sl Komaba, Tokyo 153--8914, Japan}
\\
\begin{tabular}{l}
  {\small ${}^\dagger$e-mail: {\tt sakuma@ms.u-tokyo.ac.jp}}
  \vspace{-1mm}\\
  {\small ${}^\ddagger$e-mail: {\tt yamauchi@ms.u-tokyo.ac.jp}}
\end{tabular}
}

\begin{document}

\baselineskip 6mm

\maketitle

\vspace{-5mm}

\begin{abstract}
  In this article we study and obtain a classification of Ising vectors in vertex operator algebras
  associated to binary codes and $\sqrt{2}$ times root lattices,
  where an Ising vector is a conformal vector with central charge 1/2 generating a simple
  Virasoro sub VOA.
  Then we apply our results to study certain commutant subalgebras related to root systems.
  We completely classify all Ising vectors in such commutant
  subalgebras and determine their full automorphism groups.
\end{abstract}


\section{Introduction}

Motivated by the problem of looking for maximal associative subalgebras of the monstrous
Griess algebra \cite{G}, a class of conformal vectors in the lattice
vertex operator algebra $V_{\sqrt{2}R}$ were studied and constructed in
\cite{DLMN}, where $R$ is a root lattice of type $A,D$ or $E$ of rank $\ell$
and $\sqrt{2}R$ denotes $\sqrt{2}$ times an ordinary root lattice $R$.
We adopt the standard notation for lattice vertex operator algebras as in \cite{FLM}.
In \cite{DLMN}, Dong et al.\ constructed conformal vectors of $V_{\sqrt{2}R}^+$ of
the forms
$$
  s_R=\dfr{h}{h+2}\w -\dfr{1}{h+2} \sum_{\alpha\in \Phi(R)}e^{\sqrt{2}\alpha}
  \q \text{and}\q
  \vsb\\
  \tilde{\w}_R= \dfr{2}{h+2}\w +\dfr{1}{h+2} \sum_{\alpha\in \Phi(R)} e^{\sqrt{2}\alpha},
$$
where $\w$ is the Virasoro element of $V_{\sqrt{2}R}^+$, $h$ is the
Coxeter number of $R$ and $\Phi(R)$ denotes the root system of $R$.
The central charges of $s_R$ and $\tilde{\w}_R$ are respectively
$\ell h/(h+2)$ and $2\ell/ (h+2)$. The Weyl group $W(R)$ of the root
system $\Phi (R)$ induces a natural action on the lattice vertex
operator algebra $V_{\sqrt{2}R}$ and its $\Z_2$-orbifold
$V_{\sqrt{2}R}^+$. By the construction, both conformal vectors $s_R$
and $\tilde{\om}_R$ are fixed by $W(R)$ so that $W(R)$ acts
identically on the Virasoro vertex operator subalgebra $\vir
(\tilde{\w}_R)$ generated by $\tilde{\w}_R$. Therefore, the
commutant (or coset) subalgebra
$$
  \com_{V_{\sqrt{2}R}^+}(\vir(\tilde{\w}_R))
   =\{ a\in V_{\sqrt{2}R}^+ \mid [Y(a,z_1), Y(u,z_2)]=0\ \text{for all}\
         u\in \vir (\tilde{\w}_R)\}
$$
affords an action of the Weyl group $W(R)$.

In this article, we shall study the structure of the commutant vertex operator
subalgebra $\com_{V_{\sqrt{2}R}^+}(\vir(\tilde{\w}_R))$ of $V_{\sqrt{2}R}^+$.
As our main result, we shall determine all the conformal vectors of central
charge $1/2$ in $\com_{V_{\sqrt{2}R}^+}(\vir(\tilde{\w}_R))$ and show that
the vertex operator algebra $\com_{V_{\sqrt{2}R}^+}(\vir(\tilde{\w}_R))$ is
generated by its weight two subspace.
We shall also determine the full automorphism group of
$\com_{V_{\sqrt{2}R}^+}(\vir(\tilde{\w}_R))$ and show that it always contains
a half of the Weyl group $W(R)$.

Besides its own interest, the study of
$\com_{V_{\sqrt{2}R}^+}(\vir(\tilde{\w}_R))$ also has a perspective
on the Monster simple group. There is an attempt in \cite{LYY1,LYY2}
to elucidate McKay's observation on the Monster simple group and the
extended $E_8$ diagram via vertex operator algebras. It is known
(cf.\ \cite{C,M1}) that a 2A-involution of the Monster is in
one-to-one correspondence with a conformal vector with central
charge 1/2 in the moonshine vertex operator algebra \cite{FLM} via
so-called Miyamoto involution. Taking notice of a fact that the
conformal vector $\tilde{\w}_{E_8}$ of $V_{\sqrt{2}E_8}^+$ is of
central charge 1/2, Lam et al.\ used $V_{\sqrt{2}E_8}$ to relate the
$E_8$ diagram with the Monster in \cite{LYY1,LYY2}. To each node of
the extended $E_8$ diagram, in \cite{LYY1} they constructed two
conformal vectors of central charge 1/2 inside $V_{\sqrt{2}E_8}$
which shall induce two 2A-involutions  of the Monster. Then they
used commutant subalgebra $\com_{V_{\sqrt{2}R}^+}(\vir(s_R))$ for a
sublattice $R$ of $E_8$ to study the subalgebra generated by two
such conformal vectors of $V_{\sqrt{2}E_8}$ in \cite{LYY2}. Since
$s_R$ and $\tilde{\w}_R$ are mutually commutative conformal elements
in $V_{\sqrt{2}R}^+$ and the sum $s_R+\tilde{\w}_R$ is the Virasoro
vector of $V_{\sqrt{2}R}^+$, the subalgebra
$\com_{V_{\sqrt{2}R}^+}(\vir(\tilde{\w}_R))$ is exactly equal to the
commutant subalgebra of $\com_{V_{\sqrt{2}R}^+}(\vir(s_R))$ in
$V_{\sqrt{2}R}^+$ (cf.\ \cite{FZ}). By this duality, we can relate
the structure of $\com_{V_{\sqrt{2}R}^+}(\vir(\tilde{\w}_R))$ with
the centralizer of two 2A-involutions of the Monster. In \cite{GN},
Glauberman and Norton suggested some interesting relations between
Weyl groups and centralizers of two 2A-involutions of Monster. Since
$\com_{V_{\sqrt{2}R}^+}(\vir(\tilde{\w}_R))$ naturally affords an
action of the Weyl group $W(R)$, the study of
$\com_{V_{\sqrt{2}R}^+}(\vir(\tilde{\w}_R))$ may lead
Glauberman-Norton's observation to an appreciable settlement.


Let us denote $\com_{V_{\sqrt{2}R}^+}(\vir(\tilde{\w}_R))$ by $M_R$ for simplicity of
notation.
The structure of $M_R$ is closely related to that of the root system $\Phi(R)$ of $R$,
so the automorphism group of $M_R$ has a similar structure to the Weyl group $W(R)$ of $R$.
The Weyl group $W(R)$ is a 3-transposition group, and Miyamoto discovered in \cite{M1}
that the 3-transposition property of $W(R)$ acting on $V_{\sqrt{2}R}^+$ comes from
the structure of conformal vectors with central charge 1/2.
In \cite{M1}, Miyamoto introduce a way to define involutions of vertex operator algebras
containing conformal vectors with central charge 1/2.
He showed that in some cases these involutions generate a 3-transposition group,
which is exactly the case for $V_{\sqrt{2}R}^+$.
There are many results about the group generated by these involutions, see
\cite{KM, Ma, M1, La1, LS, LYY1, LYY2, Y1}.
Recently, Matsuo classified all 3-transposition groups defined by
conformal vectors with central charge 1/2 in \cite{Ma}.
According to his classification, the parameters of $M_R$ such as central charges
or dimensions of weight two subspace coincide with those for 3-transposition groups
in his list.
In particular, if we take $R=E_6,E_7$, or $E_8$, then the corresponding groups
have nice symmetry.
In fact, this is one of the main motivations of our work.
We shall determine the automorphism group of $M_R$ by studying conformal vectors
of central charge 1/2 inside $M_R$.

Let us explain our result more precisely.
Since the conformal vector with central charge 1/2 plays a central role in our discussion,
we will refer a conformal vector with central charge 1/2 to as an {\it Ising vector} if
it generates a simple Virasoro vertex operator subalgebra.
To study Ising vectors, we first study a class of vertex operator algebras, called
{\it code vertex operator algebras} (cf.\ \cite{M2}).
In the case that a code vertex operator algebra has finitely many Ising vectors,
we classify all the Ising vectors of a code vertex operator algebra.
Moreover, we will present a formula to count the exact number of such conformal
vectors in Corollary \ref{cor:3.8}.
Then for a (not necessarily indecomposable) root lattice $R$ we will classify all
Ising vectors inside $V_{\sqrt{2}R}^+$ by giving an embedding of $V_{\sqrt{2}R}^+$ into
a code vertex operator algebra.
As a result, it is shown in Theorem \ref{thm:4.6} that for each Ising vector $e$ of
$V_{\sqrt{2}R}^+$, there exists a sublattice $K$ of $\sqrt{2}R$ isometric to either
$\sqrt{2}A_1$ or $\sqrt{2}E_8$ such that $e\in V_{K}^+\subset V_{\sqrt{2}R}^+$.
We expect that this is true not only for a lattice vertex operator algebra associated to
$\sqrt{2}R$ but also for any lattice without roots (cf.\ Remark \ref{rem:4.7}).
The classification of Ising vectors of $M_R$ immediately follows from that of
$V_{\sqrt{2}R}^+$.
We also show that $M_R$ is generated by its Ising vectors.
This reduces the study of $\aut(M_R)$ to an analysis of the permutation group of
the Ising vectors.
Since we know the group generated by involutions associated to these Ising vectors
by the results obtained in \cite{KM,Ma}, it is not difficult to determine
the permutation group of the Ising vectors.

The organization of this paper is as follows.
In Section 2 we prepare some basic notation and related facts about conformal vectors
and commutant subalgebras.
In Section 3 we study Ising vectors of code vertex operator algebras.
We present a classification of Ising vectors of a code vertex operator algebras
in Proposition \ref{prop:3.7}.
In Section 4 we study Ising vectors of lattice vertex operator algebras.
We classify Ising vectors of $V_{\sqrt{2}R}^+$ and $M_R$ with $R$ an irreducible root lattice.
In Section 5 we determine the shape of $\aut (M_R)$.
We make use of the fact that $M_R$ is generated by its Ising vectors.
In the case that $R=E_6$ or $E_7$, the proof of this fact is rather technical and long.
So we separately give it in Section 7.
In Section 6 we study a relation between inductive structures of a 3-transposition group
acting on a vertex operator algebra and its commutant subalgebras.
As an example, we study a vertex operator algebra corresponding to an inductive structure
$\mathrm{O}_{10}^+(2)^{(2)}\simeq \mathrm{O}_8^-(2)$.
Section 7 is devoted to the proof that $M_R$ with $R=E_6$ or $E_7$ is generated by
its weight two subspace as a vertex operator algebra.
We use some representation theory of the unitary Virasoro vertex operator algebras
and $W$-algebras there.
So we also review some facts about these vertex operator algebras.

\paragraph{Acknowledgement}
The third author (H.Y) thanks Hiroki Shimakura for discussions, especially for
valuable comments on Ising vectors of lattice vertex operator algebras and the
moonshine vertex operator algebra.
Part of the work was done when the second and the third authors
were visiting the National Center for Theoretical Sciences, Taiwan
on September 2004. They thank W.F.\ Ke for his hospitality.

\subsection{Notation}
In this article, every vertex operator algebra (VOA for short) we considered
is defined over $\C$ and supposed to have the grading $V=\oplus_{n\geq 0}V_n$
with $V_0=\C \vac$.
For a VOA structure $(V,Y(\cd,z),\vac,\w)$ on $V$, the vector $\w$ is called  the
{\it Virasoro element} or {\it Virasoro vector} of $V$.
For simplicity, we often use $(V,\w)$ to denote the structure $(V,Y(\cd,z),\vac,\w)$.
The vertex operator $Y(a,z)$ of $a\in V$ is expanded such as
$Y(a,z)=\sum_{n\in \Z}a_{(n)} z^{-n-1}$.

An element $u\in V$ is referred to as a {\it conformal vector} with central charge
$c_u\in \C$ if $u\in V_2$ and it satisfies $u_{(1)}u=2u$ and $u_{(3)}u=c_u\vac$.
It is well-known (cf.\ \cite{M1,La1}) that after setting $L^u(n):=u_{(n+1)}$, $n\in \Z$,
we obtain a representation of the Virasoro algebra on $V$:
$$
  [L^u(m),L^u(n)]=(m-n)L^u(m+n)+\delta_{m+n,0}\dfr{m^3-m}{12}c_u.
$$
For $c,h\in \C$, we denote by $L(c,h)$ the irreducible highest weight module over the
Virasoro algebra with central charge $c$ and highest weight $h$.
It is well-known that $L(c,0)$ has a simple VOA structure (cf.\ \cite{FZ}).

For a positive definite even lattice $L$, we will denote the lattice VOA
associated to $L$ by $V_L$ (cf.\ \cite{FLM}).
We adopt the standard notation for $V_L$ as in \cite{FLM}.
In particular, $V_L^+$ denotes the fixed point subalgebra of $V_L$ under
the lift of $(-1)$-isometry on $L$.

Given an automorphism group $G$ of $V$, we denote by $V^G$ the $G$-fixed point
subalgebra of $V$.
The subalgebra $V^G$ is called the {\it $G$-orbifold} of $V$ in the literature.

We denote the ring $\Z/2\Z$ by $\Z_2$.
Let $C\subset \Z_2^n$ be a linear code.
For a codeword $\alpha=(\alpha_1,\dots,\alpha_n)\in C$, we define the support of
$\alpha$ by $\supp(\alpha):=\{ i \mid \alpha_i=1\}$.
For a subset $A$ of $C$, we define $\supp(A):= \cup_{\alpha\in A} \supp(\alpha)$.

\section{Commutant subalgebras and conformal vectors}

We will present the notion of a commutant subalgebra and its description
by a pair of mutually commutative conformal vectors.

\subsection{Commutant subalgebras}

Let $(V,\w)$ be a VOA. A decomposition $\w=u^1+\cds +u^n$ is called
{\it orthogonal} if $u^i$ are conformal vectors and they are mutually commutative,
i.e.,  $[Y(u^i,z_1),Y(u^j,z_2)]=0$ if $i\ne j$.
The following lemma is well-known (cf.\ Theorem 5.1 of \cite{FZ}).

\begin{lem}\label{lem:2.1}
  Assume that a VOA $(V,\w)$ has a grading $V=\oplus_{n\geq 0}V_n$ such that
  $V_0=\C \vac$ and $V_1=0$.
  Then for any conformal vector $u$ of $V$, the decomposition $\w=u+(\w-u)$ is
  orthogonal.
\end{lem}

In this paper, a sub VOA of $V$ is defined by means of a pair $(U,e)$ of
a subalgebra $U$ of $V$ containing the vacuum element $\vac$ of $V$ and
a conformal vector $e$ of $V\cap U$ such that $e$ is the Virasoro vector of $U$
and the structure $(U,e)$ inherits the grading of $V$, that is,
$U=\oplus_{n\geq 0} U_n$ with $U_n=V_n\cap U$.
In the case that the conformal vector of a sub VOA of $V$ coincides with that
of $V$, we will refer such a sub VOA to as a {\it full} sub VOA.
Every conformal vector $u$ of $V$ together with the vacuum $\vac$ generates a
Virasoro sub VOA of $V$ which we will denote by $\vir(u)$.
\vsb

Let $S$ be a subset of $V$.
The following is easy to see.

\begin{lem}\label{lem:2.2}
  The subspace $S^c:=\{ a\in V\mid a_{(i)}S=0,\ i\geq 0\}$ forms a subalgebra and
  satisfies that $[Y(u,z_1),Y(v,z_2)]=0$ for any $u\in S$, $v\in S^c$.
\end{lem}

By the lemma above, we define the {\it commutant subalgebra} of a subalgebra $U$
of $V$ by
\begin{equation}\label{eq:2.1}
  \com_V(U):= U^c=\{ a\in V \mid a_{(i)} U=0,\ i\geq 0\} .
\end{equation}
If $U$ has a conformal vector $e$ such that $(U,e)$ forms a sub VOA of $V$,
then it is shown in Theorem 5.2 of \cite{FZ} that $\com_V(U)=\ker_V e_{(0)}$.
Therefore, the commutant subalgebra is described in term of the conformal  vector
$e$ of $U$.

\subsection{Conformal vectors associated to root systems}

Let $R$ be a root lattice with root system  $\Phi(R)$.
Let $\ell$ be the rank of $R$ and $h$ the Coxeter number of $R$.
We denote by $\sqrt{2} R$ the lattice whose norm is twice of $R$'s.
We consider the fixed point subalgebra $V_{\sqrt{2}R}^+$ of the lattice VOA
$V_{\sqrt{2}R}$ under the lift of $(-1)$-isometry on $R$.
It is clear that $V_{\sqrt{2}R}^+$ has a grading
$V_{\sqrt{2}R}^+=\oplus_{n\geq 0}(V_{\sqrt{2}R}^+)_n$ such that
$(V_{\sqrt{2}R}^+)_0=\C \vac$ and $(V_{\sqrt{2}R}^+)_1=0$.
We can find conformal vectors of $V_{\sqrt{2}R}^+$ defined as follows.
Set
\begin{equation}\label{eq:2.2}
  s=s_R:=\fr{h}{h+2}\w -\fr{1}{h+2}\sum_{\alpha \in \Phi(R)} e^{\sqrt{2} \alpha}
  \in V_{\sqrt{2}R}^+,
\end{equation}
where $\w$ is the Virasoro vector of $V_{\sqrt{2}R}^+$.
Then it is shown in \cite{DLMN} that $s$ defines a conformal vector with central charge
$\ell h/(h+2)$.
By Lemma \ref{lem:2.1},
\begin{equation}\label{eq:2.3}
  \tilde{\w}=\tilde{\w}_R:=\w-s=\fr{2}{h+2}\w+\fr{1}{h+2}\sum_{\alpha \in \Phi(R)}
  e^{\sqrt{2} \alpha} \in V_{\sqrt{2}R}^+
\end{equation}
is also a conformal vector  with central charge $2\ell/(h+2)$ and
the decomposition $\w=s+\tilde{\w}$ is orthogonal.
In this article, we will mainly consider the commutant subalgebra
\begin{equation}\label{eq:2.4}
  M_R:= \com_{V_{\sqrt{2}R}^+}\l(\vir(\tilde{\w})\r)
      = \ker_{V_{\sqrt{2}R}^+} \tilde{\w}_{(0)} .
\end{equation}
It is clear from the expression in \eqref{eq:2.3} that $\tilde{\w}$ is
invariant under the natural action of the Weyl group $W(R)$ associated to
the root system $\Phi(R)$.
Therefore, the commutant subalgebra $M_R$ naturally affords an action of
the Weyl group $W(R)$.

\section{Ising vectors and Ising frames}

We will introduce the notion of an Ising vector by which one can define an involution
of a vertex operator algebra.
We will review basic facts about involutions associated to Ising vectors.
Then we will study the classification of Ising vectors of a code vertex operator algebra.
We will also study automorphisms of a code vertex operator algebra.
A brief description of the autormophism group of a code vertex operator
algebra will be presented.

\subsection{Miyamoto involution}

We begin by the definition of an Ising vector.

\begin{df}
  A conformal vector  $e$ of a VOA $V$ is called  an {\it Ising vector}
  if the subalgebra $\vir(e)$ generated by $e$ is isomorphic to the simple Virasoro VOA
  $L(\shf,0)$ with central charge 1/2.
  An orthogonal decomposition $\w=e^1+\cds +e^n$ of the Virasoro vector $\w$ is called
  an {\it Ising frame} if each $e^i$ is an Ising vector.
\end{df}

\begin{rem}
  An Ising vector is often referred to as a rational conformal vector of central charge
  $1/2$ in the literature (cf.\ \cite{M1,La1}).
\end{rem}

It is well-known that the Virasoro VOA $L(\shf,0)$ is rational and has three irreducible
representations, $L(\shf,0)$, $L(\shf,\shf)$ and $L(\shf,\sfr{1}{16})$ (cf.\ \cite{DMZ}).

Let $e$ be an Ising vector of a VOA $V$.
Since $\vir(e)$ is rational, $V$ is a semisimple $\vir(e)$-module.
For $h=0,1/2,1/16$, denote by $V_e(h)$  the sum of all irreducible
$\vir(e)$-submodules of $V$ isomorphic to $L(\shf,h)$.
Then we have the isotypical decomposition:
\begin{equation}\label{eq:3.1}
  V=V_e(0)\oplus V_e(1/2)\oplus V_e(1/16).
\end{equation}
Define a linear automorphism $\tau_e$ on $V$ which acts on $V_e(0)\oplus V_e(\shf)$ by
identity and on $V_e(\sfr{1}{16})$ by $-1$.
Then it is shown in \cite{M1} that $\tau_e\in \aut (V)$.
On the $\la \tau_e\ra$-fixed point subalgebra $V^{\la \tau_e\ra}=V_e(0)\oplus V_e(\shf)$,
define a linear automorphism $\sigma_e$ which acts on $V_e(0)$ by identity and
on $V_e(\shf)$ by $-1$.
Then it is also shown in \cite{M1} that $\sigma_e\in \aut ( V^{\la \tau_e\ra})$.
We will refer $\tau_e\in \aut(V)$ (resp.\ $\sigma_e\in \aut(V^{\la \sigma_e\ra})$) to
as the {\it Miyamoto involution} of $\tau$-type (resp.\ $\sigma$-type).
An Ising vector $e$ of $V$ is called {\it of $\sigma$-type} if $\tau_e$ defines identity
on $V$, and we also refer an Ising frame $\w=e^1+\cds +e^n$ of $V$ to as {\it of $\sigma$-type}
if all $e^i$, $1\leq i\leq n$, are of $\sigma$-type on $V$.

\subsection{Code VOA}

Let us review the construction of code VOAs in \cite{M2} at least for what we need
in this paper.

Let $\mathcal{A}$ be the algebra generated by $\{ \psi_r \mid r\in \Z+1/2\}$ subject
to the defining relation $\psi_r\psi_s+\psi_s\psi_r=\delta_{r+s,0}$, $r,s\in \Z+1/2$.
Let $\mathcal{A}^+$ be the subalgebra of $\mathcal{A}$ generated by
$\{ \psi_r \mid r>0\}$ and let $\C \vac$ be a trivial $\mathcal{A}$-module.
Then set $X:=\ind_{\mathcal{A}^+}^{\mathcal{A}} \C \vac$.
Consider the generating function
\begin{equation}\label{eq:3.2}
  \psi(z):=\sum_{n\in \Z} \psi_{n+1/2}z^{-n-1} .
\end{equation}
It is well-known that the space $X$, with the standard $\Z_2$-grading, has a unique
structure of a simple vertex operator superalgebra (SVOA for short) with
the vacuum element $\vac$ such that $Y(\psi_{-1/2}\vac,z)=\psi(z)$.
The vector $\w=\hf \psi_{-3/2}\psi_{-1/2}\vac$ is a Virasoro vector of $X$ with
central charge 1/2 and the quadruple $(X,Y(\cd,z),\vac,\w)$ is isomorphic to
$L(\shf,0)\oplus L(\shf,\shf)$ (cf.\ \cite{KR}).

Set $X^0:=L(\shf,0)\subset X$ and $X^1:=L(\shf,\shf)\subset X$
under the isomorphism $X\simeq L(\shf,0)\oplus L(\shf,\shf)$.
Then $X^{\tensor n}$ also forms an SVOA as a tensor product of SVOAs.
For an even linear subcode $C$ of $\Z_2^n$, set
\begin{equation}\label{eq:3.3}
  V_C:=\bigoplus_{\alpha=(\alpha_1,\dots,\alpha_n) \in C}
  X^{\alpha_1}\tensor \cds \tensor X^{\alpha_n},
\end{equation}
which is a subalgebra of $X^{\tensor n}$.
This is a simple VOA called a {\it code VOA} associated to $C$ (cf.\ \cite{M2}).
A code VOA has the standard Ising frame of $\sigma$-type and the following theorem
characterizes code VOAs via their standard Ising frames.

\begin{thm}\label{thm:3.2}
  (\cite{M3})
  Let $V$ be a simple VOA with an Ising frame $\w=e^1+\cds +e^n$ of $\sigma$-type.
  Then there exists a unique even linear subcode $C$ of $\Z_2^n$ such that $V$ is isomorphic to
  a code VOA $V_C$ with respect to the Ising frame $\w=e^1+\cds +e^n$.
\end{thm}

We will use the following notation for code VOAs.
Set $u^0:=\vac\in X^0$, $u^1:=\psi_{-1/2}\vac\in X^1$, and for a codeword
$\alpha=(\alpha_1,\dots,\alpha_n)\in \Z_2^n$, we set
\begin{equation}\label{eq:3.4}
  x^\alpha:=u^{\alpha_1}\tensor \cds \tensor u^{\alpha_n} \in X^{\tensor n}.
\end{equation}
Then $x^\alpha$ is a highest weight vector of $X^{\alpha_1}\tensor \cds \tensor X^{\alpha_n}$
with norm $\la x^\alpha,x^\alpha\ra=\pm 1$\footnote{
The signs depend on the choice of 2-cocycle which we have use to construct a tensor product
$X^{\tensor n}$.},
where the invariant bilinear form $\la\cd,\cd\ra$ on $X$ is normalized such as
$\la \vac,\vac\ra =1$.

For a codeword $\beta\in \Z_2^n$, the subspace
$$
  V_{C+\beta}:=\bigoplus_{\gamma=(\gamma_1,\dots,\gamma_n)\in C+\beta}
  X^{\gamma_1}\tensor \cds \tensor X^{\gamma_n}
$$
uniquely forms an irreducible $V_C$-submodule of $X^{\tensor n}$ (cf.\ \cite{M3}).
We will call $V_{C+\beta}$ a {\it coset type module} over $V_C$.

\subsection{Ising vectors of $\sigma$-type}

In this subsection we consider a VOA $V$ with trivial weight one subspace, i.e.,  $V_1=0$.
Then the weight two subspace $V_2$ equipped with the product $a\cd b:=a_{(1)}b$ for $a,b\in V_2$
forms a commutative algebra with an invariant bilinear form $\la\cd,\cd\ra$ defined by
$\la a,b\ra\vac=a_{(3)}b$.
This algebra is called the {\it Griess algebra} of $V$.
Let $e,f$ be Ising vectors of $V$ of $\sigma$-type.
Then $e/2$ and $f/2$ are idempotents of the Griess algebra $V_2$ with squared norm 1/16.
The following result is fundamental:

\begin{prop}\label{prop:3.3}
  (\cite{MM}\cite{M1})
  Let $e,f$ be distinct Ising vectors of $V$ of $\sigma$-type.
  Then one of the following holds:
  \\
  (1) $\la e,f\ra= 0$ and $e\cd f=0$.
  In this case $\sigma_e f=f$ and $\sigma_f e=e$.
  \\
  (2) $\la e,f\ra=1/32$ and $e\ne \sigma_f e=\sigma_e f\ne f$.
  In this case $\sigma_e\sigma_f$ is of order 3 and the following equality holds in $V_2$:
  \begin{equation}\label{eq:3.5}
    e\cd f =\fr{1}{4}(e+f-\sigma_e f).
  \end{equation}
  Therefore, $e$ and $f$ generates a three dimensional subalgebra $\C e\oplus \C f\oplus
  \C \sigma_ef $ in the Griess algebra $V_2$ on which the symmetric group of degree three
  acts.
\end{prop}

Let $H_8$ be the [8,4,4]-Hamming code:
$$
  H_8:=\Span_{\Z_2}\{ (1111 1111), (11110000), (11001100), (10101010) \} \subset \Z_2^8.
$$
It is well-known that $H_8$ is the unique doubly even self-dual code of length 8.
Let $V_{H_8}$ be the code VOA associated to $H_8$ and
let $\w=e^1+\cds +e^8$ be the standard Ising frame of $V_{H_8}$.
Inside the Hamming code VOA $V_{H_8}$, we can find three Ising frames.
For a codeword $\alpha \in \Z_2^8$, define
\begin{equation}\label{eq:3.6}
  t^\alpha:= \dfr{1}{8}\sum_{i=1}^8 e^i +
  \dfr{1}{8}\sum_{\beta\in H_8,\, \la \beta,\beta\ra=4} (-1)^{\la \alpha,\beta\ra} x^\beta
  \in V_{H_8},
\end{equation}
where $x^\beta$ above are defined as in \eqref{eq:3.4}.
Then it is shown in \cite{M2} that $t^\alpha$ is an Ising vector of $\sigma$-type.

\begin{prop}\label{prop:3.4}
  (\cite{MM}\cite{M4})
  Inside the Hamming code VOA $V_{H_8}$, there are exactly three Ising frames given as follows.
  \begin{equation}\label{eq:3.7}
    I_0:=\{ e^i \mid 1\leq i\leq 8\},\
    I_1:=\{ t^{\nu^j} \mid 1\leq j\leq 8\}\ \mathrm{and }\
    I_2:=\{ t^{\nu^1+\nu^k} \mid 1\leq k\leq 8\},
  \end{equation}
   where we have set $\nu^1:=(1000 0000)$,
   $\nu^2:=(0100 0000),\dots,\nu^8:=(0000 0001)\in \Z_2^8$.
   Moreover, if $f\in I_a$, then $\sigma_f I_b=I_c$ if $\{a,b,c\}=\{0,1,2\}$
   so that all the frames are mutually conjugate to each others.
\end{prop}

Based on Propositions \ref{prop:3.3} and \ref{prop:3.4}, the following result is
established in \cite{La1}.

\begin{prop}\label{prop:3.5}
  (\cite{La1})
  Let $C$ be an even linear code whose minimum weight is greater than 2.
  If the code VOA $V_C$ contains an Ising vector $f$ of $\sigma$-type which is
  not a summand of the standard Ising frame of $V_C$, then $C$ contains a subcode
  $D$ isomorphic to $H_8$ such that $f\in V_D\subset V_C$ and $f$ is of the form
  \eqref{eq:3.6} in $V_D$.
  In particular, there are exactly 24 Ising vectors inside $V_{H_8}$ and every
  Ising vector of $V_{H_8}$ is a summand of an Ising frame of $V_{H_8}$.
\end{prop}

We give a generalization of the proposition above.
The following lemma enables us to reduce a general case to the known case.

\begin{lem}\label{lem:3.6}
  Let $V$ be a VOA with $V_1=0$.
  Suppose that $V$ has two Ising vectors $e,f$ and $e$ is of $\sigma$-type.
  Then $e\in V^{\la \tau_f\ra}$.
\end{lem}

\pf
Consider the Griess algebra $V_2$ of $V$.
Let $t$ be an Ising vector of $V$.
Then it is shown in \cite{M1} that the Griess algebra $V_2$ affords the
orthogonal decomposition
\begin{equation}\label{eq:3.8}
  V_2= \C t\perp (V_t(0)\cap V_2) \perp (V_t(\shf)\cap V_2)
       \perp (V_t(\sfr{1}{16})\cap V_2)
\end{equation}
and $L^t(0)=t_{(1)}$ acts on $V_t(h)\cap V_2$ by the scalar $h$.
For convention of notation, we set $B_t(h):=V_t(h)\cap V_2$ for $h=0,\shf,\sfr{1}{16}$,
and for $x\in V_2$, we will denote the corresponding decomposition by
$x=\lambda_{x,t} t +x_t(0)+x_t(\shf)+x_t(\sfr{1}{16})$ with $\lambda_{x,t}\in \C$.
Since $\la x,t\ra=\lambda_{x,t} \la t,t\ra =\lambda_{x,t}/4$, the
scalar $\lambda_{x,t}$ is given by $4\la x,t\ra$.
By the fusion rules of $L(\shf,0)$-modules (cf.\ \cite{DMZ}), $V_2$ has the
following structure:
\begin{equation}\label{eq:3.9}
\begin{array}{ll}
  B_t(0)\cd B_t(h)\subset B_t(h),\  h=0,\shf,\sfr{1}{16},\q
  & B_t(\shf)\cd B_t(\shf) \subset \C t\oplus B_e(0),
  \vsb\\
  B_t(\shf)\cd B_t(\sfr{1}{16})\subset B_t(\sfr{1}{16}),
  & B_t(\sfr{1}{16})\cd B_t(\sfr{1}{16})\subset \C t\oplus B_t(0)\oplus B_t(\shf).
\end{array}
\end{equation}
Write $f=\lambda e+f_e(0)+f_e(\shf)$ and $e=\lambda f+e_f(0)+e_f(\shf)+e_f(\sfr{1}{16})$
with $\lambda=4\la e,f\ra$.
We shall show that $e_f(\sfr{1}{16})=0$.
Since $f$ is an Ising vector, we have $f\cd f=2f$.
Using \eqref{eq:3.9}, we compare the $V_e(\shf)$-parts of both sides of $f\cd f=2f$ and
obtain
\begin{equation}\label{eq:3.10}
  f_e(0)\cd f_e(\shf)
  = \hf(2-\lambda)  f_e(\shf).
\end{equation}
Similarly, we compare the $V_f(\sfr{1}{16})$-parts in the equality $e\cd e=2e$ and
get
\begin{equation}\label{eq:3.11}
  e_f(0)\cd e_f(\sfr{1}{16})+e_f(\shf)\cd e_f(\sfr{1}{16})
  = \fr{1}{16}(16-\lambda) e_f(\sfr{1}{16}).
\end{equation}
Since the Griess algebra $V_2$ is commutative, $e\cd f=f\cd e$ and we have
\begin{equation}\label{eq:3.12}
  f_e(\shf)=4\lambda(1-\lambda)f-4\lambda e_f(0)+(1-4\lambda)e_f(\shf)
  +\dfr{1}{8}(1-32\lambda)e_f(\sfr{1}{16}).
\end{equation}
Using $f_e(0)=f-\lambda e-f_e(\shf)$, we have
\begin{equation}\label{eq:3.13}
  f_e(0)=(3\lambda-1)(\lambda-1)f+3\lambda e_f(0)+(3\lambda-1)e_f(\shf)
  +\fr{1}{8}(24\lambda-1)e_f(\sfr{1}{16}).
\end{equation}
From the $V_f(\sfr{1}{16})$-part in the equality $e\cd f_e(0)=0$ we obtain
\begin{equation}\label{eq:3.14}
  e_f(\shf)\cd e_f(\sfr{1}{16})=\fr{1}{16}(92\lambda-1)e_f(\sfr{1}{16}).
\end{equation}
Then by \eqref{eq:3.11} we have
\begin{equation}\label{eq:3.15}
  e_f(0)\cd e_f(\sfr{1}{16})=\fr{1}{16}(-93\lambda+17)e_f(\sfr{1}{16}).
\end{equation}
By \eqref{eq:3.14}, \eqref{eq:3.15} and the $V_f(\sfr{1}{16})$-parts
in the equality \eqref{eq:3.10}, we obtain
\begin{equation}\label{eq:3.16}
  \fr{1}{128}(32\lambda-1)(64\lambda+13)e_f(\sfr{1}{16})=0.
\end{equation}
Therefore, either $\lambda=1/32$, $\lambda=-13/64$ or $e_f(\sfr{1}{16})=0$.

Now suppose that $e_f(\sfr{1}{16})\ne 0$.
Then $\lambda=1/32$ or $-13/64$.
By comparing the $V_f(\shf)$-parts in the equality $e\cd e=2e$ we obtain
\begin{equation}\label{eq:3.17}
  \big( e_f(\sfr{1}{16})\cd e_f(\sfr{1}{16}) \big)_{\! f} (\shf)
  = (2-\lambda)e_f(\shf)-2e_f(0)\cd e_f(\shf).
\end{equation}
Then from the $V_f(\shf)$-part in the equality $e\cd f_e(0)=0$ together with
\eqref{eq:3.17} we get
\begin{equation}\label{eq:3.18}
  e_f(0)\cd e_f(\shf)=\fr{1}{6}(29\lambda+2)e_f(\shf)
\end{equation}
and from the $V_f(\shf)$-parts in the equality \eqref{eq:3.10} we also have
\begin{equation}\label{eq:3.19}
  \fr{1}{32}(168\lambda+1)e_f(0)\cd e_f(\shf)
  = \fr{1}{64}(600\lambda^2-17\lambda+34)e_f(\shf).
\end{equation}
Since $\lambda=1/32$ or $-13/64$, the equalities \eqref{eq:3.18} and \eqref{eq:3.19}
imply that $e_f(\shf)=0$.
Then the equation \eqref{eq:3.14} contradicts to our assumption.
Thus $e_f(\sfr{1}{16})=0$ and the lemma follows.
\qed
\vsb

By the lemma above, we can generalize Proposition \ref{prop:3.5} as follows.

\begin{prop}\label{prop:3.7}
  Let $C$ be an even linear code with minimum weight greater than 2, and let $f$ be an
  Ising vector of the associated code VOA $V_C$.
  If $f$ is not a summand of the standard Ising frame of $V_C$, then there is a subcode $D$ of
  $C$ isomorphic to $H_8$ such that $f\in V_D\subset V_C$ and $f$ is of the form \eqref{eq:3.6}
  in $V_D$.
\end{prop}

\pf
Let $\w=e^1+\cds +e^n$ be the standard Ising frame of $V_C$.
Then by Lemma \ref{lem:3.6} all $e^i$ are contained in the $\tau_f$-fixed point subalgebra
$V^{\la \tau_f\ra}_C$ of $V_C$.
Therefore, there is a subcode $C'$ of $C$ with index at most two such that
$V_C^{\la \tau_f\ra}$ is a code VOA $V_{C'}$ with respect to the Ising
frame $\w=e^1+\cds +e^n$.
Since $f$ is of $\sigma$-type on $V_{C'}$, we can apply Proposition \ref{prop:3.5} to
$V_{C'}$.
This completes the proof.
\qed

\begin{cor}\label{cor:3.8}
  Let $C$ be an even linear code of length $n$ and assume that the minimum weight of $C$ is
  greater than 2.
  \\
  (1) Every Ising vector of $V_C$ is a summand of an Ising frame of $V_C$.
  \\
  (2) Let $N$ be the number of embeddings of $H_8$ into $C$.
  Then $V_C$ contains exactly $16N+n$ Ising vectors.
\end{cor}

As an application of Proposition \ref{prop:3.5}, we show the conjugacy property of
Ising frames of a code VOA.

\begin{lem}\label{lem:3.9}
  Let $C$ be an even linear code of length $n$ with no weight two element.
  Suppose that there is a subcode $D$ of $C$ which is isomorphic to
  the Hamming code $H_8$.
  Then an Ising vector $t\in V_D\subset V_C$ of the form \eqref{eq:3.6} is
  of $\sigma$-type if and only if $\abs{\supp (\alpha)\cap \supp(D)}$ is even for
  all codeword $\alpha \in C$.
\end{lem}

\pf
We may assume that $\supp(D)=\{1,\dots,8\}$.
Let $\alpha =(\alpha_1,\dots,\alpha_n)$ be a codeword of $C$.
Then $V_C$ is a sum of a coset type $V_D$-submodules $V_{D+(\alpha_1,\dots,\alpha_8)}$
as a $V_D$-module.
It is shown in Theorem 2.2 of \cite{M4} that $\tau_t$ is trivial on
$V_{D+(\alpha_1,\dots,\alpha_8)}$ if and only if $(\alpha_1,\dots,\alpha_8)$ is an
even codeword.
So the assertion follows.
\qed
\vsb

By Proposition \ref{prop:3.5} and the lemma above, we can in principal count the
number of $\sigma$-type Ising frames of $V_C$.
Denote by $G$ the subgroup of $\aut (V_C)$ generated by $\sigma$-type Miyamoto involutions.
By Proposition \ref{prop:3.3}, $G$ is a 3-transposition group (see Definition \ref{df:5.1}).

\begin{prop}\label{prop:3.10}
  Let $C$ be an even linear code of length $n$ whose minimum weight is greater than 2 and
  let $\w=e^1+\cds +e^n$ and $\w=f^1+\cds +f^n$ be $\sigma$-type
  Ising frames of the code VOA $V_C$.
  Then there is an element $\rho \in G$ such that
  $\{ \rho e^1,\dots, \rho e^n\} = \{ f^1,\dots,f^n\}$.
\end{prop}

\pf It is enough to show the assertion in the case that $\{
e^1,\dots,e^n\}$ is the standard Ising frame of $V_C$. Set $I=\{
e^1,\dots,e^n\}$ and $J=\{ f^1,\dots,f^n\}$. We shall prove the
assertion inductively. Assume that there is an element $\rho_i \in
G$ such that $\{ f^1,\dots,f^i\} \subset \rho_i I$, where in the
case of $i=0$ we set $\rho_0=\id$. Since $I$ and $\rho_i I$ are
conjugate, the associated binary code of $V_C$ with respect to the
frame $\rho_i I$ is still isomorphic to $C$. Therefore, by
replacing $I$ by $\rho_i I$, we may assume that $\{
f^1,\dots,f^i\} \subset I\cap J$. If $f^{i+1}\not\in I$, then
there is a subcode $D$ of $C$ isomorphic to the Hamming code $H_8$
such that $f^{i+1}\in V_D\subset V_C$ and $f^{i+1}$ is of the form
\eqref{eq:3.6} in $V_D$ by Proposition \ref{prop:3.5}. Let
$\supp(D)=\{ j_1,\dots,j_8\}$. Then $\{
\sigma_{e^{j_1}}f^{i+1},\dots,\sigma_{e^{j_8}}f^{i+1}\} = \{
\sigma_{f^{i+1}}e^{j_1},\dots,\sigma_{f^{i+1}}e^{j_8}\}$ by
Proposition \ref{prop:3.4}. Therefore, $f^{i+1}\in
\sigma_{e^{j_1}}\sigma_{f^{i+1}}I$. It follows from \eqref{eq:3.6}
that $\la f^{i+1},e^j\ra=1/32$ for $j\in \supp(D)$. Hence $\{
f^1,\dots,f^i\} \cap \{ e^{j_1},\dots,e^{j_8}\}=\emptyset$ and
$\sigma_{e^{j_1}}\sigma_{f^{i+1}} \{ f^1,\dots,f^i\}=\{
f^1,\dots,f^i\}$. Thus
$\rho_{i+1}=\sigma_{e^{i_1}}\sigma_{f^{i+1}}\in G$ satisfies $\{
f^1,\dots,f^{i+1}\} \subset \rho_{i+1} I$. By this procedure, we
will obtain $\rho=\rho_n\in G$ such that $\rho I=J$. \qed

\begin{rem}
  We can slightly generalize Proposition \ref{prop:3.10} as follows.
  Let $\w=e^1+\cds +e^n$ be a $\sigma$-type Ising frame of $V_C$ and $\w=f^1+\cds +f^n$
  {\it any} Ising frame of $V_C$.
  Let $H$ be the subgroup of $\aut (V_C)$ generated by $\tau$-type Miyamoto involutions
  $\{ \tau_{f^i} \mid 1\leq i\leq n\}$.
  Then $f^i$ are $\sigma$-type Ising vectors on $V_C^H$.
  By Lemma \ref{lem:3.6}, $\{ e^1,\dots,e^n\}$ is contained in $V_C^H$.
  Therefore, by Proposition \ref{prop:3.10}, we can find $\rho \in \aut (V_C^H)$
  such that $\rho \{ e^1,\dots,e^n\} =\{ f^1,\dots,f^n\}$.
  However, we cannot find such $\rho$ inside $\aut (V_C)$ unless $H =1$.
  In fact, conjugating VOA structures by $\rho$, we can perform a $\Z_2$-twisted
  orbifold construction of a framed VOA (cf.\ \cite{M5, Y1}).
\end{rem}

We have only considered code VOAs associated to codes without weight two elements.
If an even linear code $C$ contains a weight two element, then the code VOA $V_C$ contains
a subalgebra isomorphic to a lattice VOA $V_{\Z \alpha}$ associated to
a lattice $\Z\alpha$ with $\la \alpha,\alpha\ra =4$ (cf.\ \cite{DMZ}).
In this case, we can define continuous automorphisms on $V_C$ by exponential
and hence $\aut (V_C)$ is always an infinite group.
Conversely, if $C$ contains no weight two element, then $\aut (V_C)$ is finite.
This result is established in \cite{M2} in the case that $V_C$ is considered over $\R$.
By Proposition \ref{prop:3.5}, we know that the set of $\sigma$-type Ising vectors of $V_C$
is finite so that $\aut (V_C)$ is still finite even if $V_C$ is defined over $\C$.
We give a brief description of $\aut (V_C)$ as follows.

\begin{prop}\label{prop:3.12}
  Let $C$ be an even linear code of length $n$ without weight two codewords.
  Then the automorphism group $\aut (V_C)$ of the code VOA $V_C$ is a finite group
  generated by the lift of $\aut (C)$ and $\sigma$-type Miyamoto involutions.
\end{prop}

\pf
Let $I=\{ e^1,\dots,e^n\}$ be the standard Ising frame of $V_C$ and
let $G$ be the 3-transposition subgroup of $\aut (V_C)$ generated by $\sigma$-type
Miyamoto involutions.
Take any $\phi \in \aut (V_C)$.
By Proposition \ref{prop:3.10}, there is an element $\rho \in G$ such that
$\rho \phi I =I$.
Then $\rho \phi$ defines an automorphism of $\vir(e^1)\tensor \cds \tensor \vir(e^n)$.
Therefore, $\rho\phi$ preserves the set $\{ \pm x^\alpha \in V_C \mid \alpha\in C\}$ of
normed highest weight vectors and thus there is a lift $\tilde{g}\in \aut (V_C)$ of
$g\in \aut (C)$ such that $\tilde{g} \rho \phi e^i=e^i$ for $i=1,\dots,n$.
Then by Schur's lemma, $\tilde{g} \rho \phi$ is written as a product of $\sigma_{e^i}$,
$1\leq i\leq n$.
Thus $\aut (V_C)$ is generated by $G$ and the lift of $\aut (C)$ on $V_C$.
Since there are finitely many Ising frames inside $V_C$, the argument above also shows
that $\aut (V_C)$ is finite.
\qed

Let $E$ be the set of $\sigma$-type Ising vectors of $V$.
We have defined a map $\sigma:E \to \aut(V)$ by associating the $\sigma$-type Miyamoto
involution $\sigma_e\in \aut(V)$ to each $e\in E$.
The following injectivity is shown in Lemma 2.5.2 of \cite{Ma}.

\begin{lem}\label{lem:3.13}(\cite{Ma})
  Assume that for each $e\in E$, there is an Ising vector $g\in V$ such that
  $\la e,g\ra=1/32$.
  Then $\sigma:E\to \aut(V)$ is injective.
\end{lem}

\pf
Suppose $\sigma_e=\sigma_f$ with $e,f\in E$.
By the assumption, there exists an Ising vector $g\in V$ such that $\la e,g\ra=1/32$.
It follows from Lemma \ref{lem:3.6} that  $e,f\in V^{\la \tau_g\ra}$.
Then by (2) of Proposition \ref{prop:3.3}, $g\ne \sigma_e g=\sigma_g e\ne e$
and hence $g\ne \sigma_g e=\sigma_e g=\sigma_f g$.
Thus $\la f,g\ra =1/32$ and $\sigma_f g=\sigma_g f$ again by Proposition \ref{prop:3.3}.
Hence $\sigma_g e=\sigma_e g=\sigma_f g=\sigma_g f$ showing $e=f$.
\qed

\section{Ising vectors of lattice VOAs}

We will classify Ising vectors of the $\Z_2$-orbifold
$V_{\sqrt{2}R}^+$ of a lattice $V_{\sqrt{2}R}$ associated with a
root lattice $R$. This classification immediately leads to a
classification of Ising vectors in $M_R$.

\subsection{A lattice VOA and its code}\label{sec:4.1}

Let $\Z^n$ be the standard lattice and let $\rho: \Z^n \to \Z_2^n$ be the reduction mod 2,
which is a group homomorphism.
For an even linear code $C\subset \Z_2^n$, the preimage $L_C:=\rho^{-1}(C)\subset \Z^n$ define
a sublattice of $\Z^n$.
Set $x^1:=(1,0,\dots,0)$, $x^2:=(0,1,\dots,0),\dots$, $x^n:=(0,0,\dots,1)\in \Z^n$.
Then $\{ x^1,\dots,x^n\}$ is the standard basis of $\Z^n$.
By definition $L_C$ contains a 4-frame $2\Z x^1\perp \cds \perp 2\Z x^n$
so that $V_{L_C}$ has a full sub VOA isomorphic to
$V_{2\Z x^1}\tensor \cds \tensor V_{2\Z x^n}$.
In the lattice VOA $V_{2\Z x^i}$, set
\begin{equation}\label{eq:4.1}
  w^{i\pm}:= \fr{1}{4}(x^i_{(-1)})^2\vac \pm \fr{1}{4}(e^{2x^i}+e^{-2x^i}).
\end{equation}
Then $w^{i\pm}$ are mutually orthogonal Ising vectors.
Let $\w$ be the Virasoro vector of $V_{L_{C}}$.
The lattice VOA $V_{L_C}$ contains an Ising frame
\begin{equation}\label{eq:4.2}
  \w = (w^{1-}+w^{1+})+(w^{2-}+w^{2+})+\cds +(w^{n-}+w^{n+}).
\end{equation}
It is easy to see that all $w^{i\pm}$ are of $\sigma$-type on $V_{L_C}$ so that
$V_{L_C}$ is a code VOA with respect to the frame.
It is obvious that $w^{i\pm}$ are contained in $V_{2\Z x^i}^+$ so that they are also contained
in $V_{L_C}^+$.
We define
$$
\begin{array}{l}
  D_0(C):= \{ (u_1,u_1,u_2,u_2,\dots,u_n,u_n)\in \Z_2^{2n} \mid (u_1,u_2,\dots,u_n)\in C^\perp\}
  \subset \Z_2^{2n},
  \vsb\\
  D_1(C):= D_0(C)\cup \l( D_0(C)+\gamma\r) \subset \Z_2^{2n},\q \gamma=(1010\dots 10)\in \Z_2^{2n}.
\end{array}
$$
As $\vir(w^{i-})\tensor \vir(w^{i+})$-modules, we have the following isomorphisms:
$$
\begin{array}{ll}
  V_{2\Z x^i}^+\simeq L(\shf,0)\tensor L(\shf,0),
  & V_{2\Z x^i}^-\simeq L(\shf,\shf)\tensor L(\shf,\shf),
  \vsb\\
  V_{2\Z x^i+x^i}^+\simeq L(\shf,0)\tensor L(\shf,\shf),
  & V_{2\Z x^i+x^i}^-\simeq L(\shf,\shf)\tensor L(\shf,0) .
\end{array}
$$
By the isomorphisms above, the structure codes of $V_{L_C}$ and $V_{L_C}^+$ with respect
to the frame \eqref{eq:4.2} are described as follows.

\begin{prop}\label{prop:4.1}
  (1) The lattice VOA $V_{L_C}$ with the Ising frame \eqref{eq:4.2} is isomorphic to
  the code VOA associated to $D_0(C)^\perp$.
  \\
  (2) The VOA $V_{L_C}^+$ with the Ising frame \eqref{eq:4.2} is isomorphic to
  the code VOA associated to $D_1(C)^\perp$.
\end{prop}

If we take
$$
  C_n:=\Span_{\Z_2}\{ (1^20^{2n-2}),(0^21^20^{2n-4}),\dots,(0^{2n-2}1^2),(1010\dots 10)\}
  \subset \Z_2^{2n},
$$
then we obtain $L_{C_n^\perp}\simeq \sqrt{2}D_{2n}$.
Therefore, the structure code of $V_{\sqrt{2}D_{2n}}$ as a code VOA is given
by the dual of the following code:
\begin{equation}\label{eq:c_n}
  \mathcal{C}_n
  :=\Span_{\Z_2}\{ (1^40^{4n-4}), (0^41^40^{4n-8}),\dots,(0^{4n-4}1^4), (1^20^2\cds 1^20^2),\}
  \subset \Z_2^{4n},
\end{equation}
and that of $V_{\sqrt{2}D_{2n}}^+$ is given by the dual code of the code
generated by $\mathcal{C}_n$ and $(1010\dots10)\in \Z_2^{4n}$.

If we take the [8,4,4]-Hamming code $H_8$, then $H_8^\perp=H_8$ and we obtain
$L_{H_8}\simeq \sqrt{2}E_8$.
The structure code of $V_{\sqrt{2}E_8}^+$ as a code VOA is given by $D_1(H_8)^\perp$,
where $D_1(H_8)$ has the following presentation:
$$
  D_1(H_8)=\Span_{\Z_2}\{ (1^{16}), (1^8 0^8), (1^4 0^4 1^4 0^4),
  (1^2 0^2  1^2 0^2 1^2 0^2  1^2 0^2), (\{ 10\}^8) \} \subset \Z_2^{16}.
$$
We note that the code $D_1(H_8)$ is the same as the first order Reed-Muller code
$\mathrm{RM}(1,4)$ of length $2^4$ so that the code $D_1(H_8)^\perp$ is equal to
the second order Reed-Muller code $\mathrm{RM}(2,4)$ (cf.\ \cite{CS}).

\subsection{Ising vectors of $V_{\sqrt{2}R}^+$}

Let $R$ be an indecomposable root lattice with root system $\Phi(R)$.
We give a classification of Ising vectors of $V_{\sqrt{2}R}^+$.
As in \eqref{eq:4.1}, for $\alpha \in \Phi(R)$ we set
\begin{equation}\label{eq:4.3}
  w^\pm(\alpha):= \fr{1}{8}\alpha_{(-1)}^2\vac
    \pm \fr{1}{4}\l( e^{\sqrt{2}\alpha}+e^{-\sqrt{2}\alpha}\r)
    \in V_{\sqrt{2}R}^+.
\end{equation}
It is shown in \cite{DMZ} \cite{DLMN} that $w^\pm(\alpha)$, $\alpha\in \Phi(R)$, are Ising
vectors of $\sigma$-type.

\begin{prop}\label{prop:4.2}
  Let $R$ be an indecomposable root lattice whose root system is of type $A_n$ or $D_n$.
  Then the set $\{ w^\pm (\alpha) \mid \alpha \in \Phi(R)\}$ exhausts all the Ising vectors
  of $V_{\sqrt{2}R}^+$.
\end{prop}

\pf
It suffices to show the assertion in the case of $R=D_{2n}$ since
$A_k$ and $D_{2m+1}$, $1\leq k,m\leq n-1$, are sublattices of $D_{2n}$.
Let $\mathcal{C}$ be the dual code of the following code:
$$
  \Span_{\Z_2}\{ (1^40^{4n-4}), (0^41^40^{4n-8}),\dots,(0^{4n-4}1^4), (1^20^2\cds 1^20^2),
  (1010\cds 1010)\} \subset \Z_2^{4n}.
$$
As shown in Sec.\ref{sec:4.1}, $V_{\sqrt{2} D_{2n}}^+$ is isomorphic to a code VOA
associated to $\mathcal{C}$.
By Corollary \ref{cor:3.8}, it suffices to count the number $N$ of embeddings of $H_8$ into
$\mathcal{C}$ and one can easily show that $N=\binom{n}{2}=n(n-1)/2$.
Therefore, there are $4n(2n-1)$ Ising vectors inside $V_{\sqrt{2}D_{2n}}^+$.
Since there are $\abs{\Phi(D_{2n})}=4n(2n-1)$ Ising vectors of the form \eqref{eq:4.3} in
$V_{\sqrt{2}D_{2n}}^+$, they are all.
\qed
\vsb

Next we consider the $E$-series.
Since $E_6$ and $E_7$ are sublattices of $E_8$, we only need to consider $V_{\sqrt{2}E_8}^+$.
Since the number of roots in $\Phi(E_8)$ is 240, we have 240 Ising vectors of $V_{\sqrt{2}E_8}^+$
of the form \eqref{eq:4.3}.
We also have an Ising vector $\tilde{\w}$ of $V_{\sqrt{2}E_8}^+$ defined as in \eqref{eq:2.3}.
Since $\tau_{\tilde{\w}}=\id$ on $V_{\sqrt{2}E_8}^+$, $\tilde{\w}$ is of $\sigma$-type.
For $h\in \C E_8$, set
\begin{equation}\label{eq:4.4}
  \varphi_h:=\exp\l( \fr{\pi\sqrt{-2}}{2} (h_{(-1)}\vac)_{(0)}\r) \in \aut (V_{\sqrt{2}E_8}).
\end{equation}
Then $\varphi_h V_{\sqrt{2}E_8}^+=V_{\sqrt{2}E_8}^+$ if and only if $\varphi_{2h}=1$.
Therefore, we obtain a group homomorphism $\varphi :E_8\ni x\mapsto \varphi_x\in
\aut(V_{\sqrt{2}E_8}^+)$ with $\ker \varphi =2E_8$.
So we have $2^8=256$ Ising vectors of $V_{\sqrt{2}E_8}^+$ of the form
\begin{equation}\label{eq:4.5}
  \varphi_x \tilde{\w},\q x\in E_8.
\end{equation}
It is shown in \cite{CS} that $E_8/2E_8$ contains 1 class represented by $0$,
120 classes represented by a pair of roots $\pm \alpha\in \Phi(E_8)$, and
135 classes represented by 16 vectors forming a 4-frame of $E_8$.

It is shown by Shimakura \cite{S2} that these $240+256=496$
vectors exhaust all the Ising vectors of $V_{\sqrt{2}E_8}^+$.
We give another proof here.

\begin{prop}\label{prop:4.3}
  There are 496 Ising vectors inside $V_{\sqrt{2}E_8}^+$.
\end{prop}

\pf
As shown in Sec.\ \ref{sec:4.1}, there is an Ising frame of $V_{\sqrt{2}E_8}^+$ such that
$V_{\sqrt{2}E_8}^+$ is  isomorphic to a code VOA associated to the Reed-Muller code
$\mathrm{RM}(2,4)$ with respect to the frame.
The minimum weight of $\mathrm{RM}(2,4)$ is 4 (cf.\ \cite{CS}).
So by Corollary \ref{cor:3.8}, it suffices to count embeddings of
the [8,4,4]-Hamming code $H_8$ to $\mathrm{RM}(2,4)$.
It is not difficult to see that for each embedding $H_8\hookrightarrow \mathrm{RM}(2,4)$,
the support of $H_8$ belongs to $\mathrm{RM}(2,4)^\perp=\mathrm{RM}(1,4)$ so that
there are exactly 30 embeddings of $H_8$ into $\mathrm{RM}(2,4)$.
Hence, there exists $30\times 16=480$ Ising vectors of the form \eqref{eq:3.6} inside
$V_{\sqrt{2}E_8}^+$.
Adding the 16 Ising vectors of summands of the standard Ising frame, there are $480+16=496$
Ising vectors which exhaust all the Ising vectors inside $V_{\sqrt{2}E_8}^+$.
We also note that all the Ising vectors are of $\sigma$-type on $V_{\sqrt{2}E_8}^+$.
\qed
\vsb

Combining Proposition \ref{prop:4.2} and \ref{prop:4.3}, we obtain the following.

\begin{thm}\label{thm:4.4}
  Let $R$ be an indecomposable root lattice with root system $\Phi(R)$.
  \\
  (1) If $R\ne E_8$, the set $\{ w^\pm (\alpha) \mid \alpha \in \Phi(R)\}$ exhausts all
  the Ising vectors of $V_{\sqrt{2}R}^+$.
  \\
  (2) There are 256 Ising vectors of $V_{\sqrt{2}E_8}^+$
  other than $\{ w^\pm (\alpha) \mid \alpha \in \Phi(E_8)\}$.
  The set $\{ w^\pm (\alpha), \varphi_x\tilde{\w} \mid \alpha\in \Phi(E_8), x\in E_8/2E_8\}$
  exhausts all the Ising vectors of $V_{\sqrt{2}E_8}^+$.
\end{thm}

\begin{cor}\label{cor:4.5}
  There is no $\tau$-type Ising vector inside the $\Z_2$-orbifold subalgebra $V_{\sqrt{2}R}^+$
  of the lattice VOA $V_{\sqrt{2}R}$ associated to an indecomposable root lattice $R$.
\end{cor}

We have treated only indecomposable root lattices.
Now we classify all Ising vectors of the $\Z_2$-orbifold $V_{\sqrt{2}L}^+$ of
a lattice VOA $V_{\sqrt{2}L}$ associated to {\it any} root lattice $L$.

\begin{thm}\label{thm:4.6}
  Let $L$ be a root lattice.
  If $e$ is an Ising vector of $V_{\sqrt{2}L}^+$, then there exists
  a sublattice $K$ of $L$ isometric to $A_1$ or $E_8$ such that
  $e\in V_{\sqrt{2}K}^+\subset V_{\sqrt{2}L}^+$.
\end{thm}

\pf
Since $L$ is a root lattice, $L$ is a direct sum of indecomposable
root lattices of {\it ADE}-type.
Let $L^{(1)}$ be the sum of irreducible sublattices of $L$ of {\it AD}-type
and $L^{(2)}$ the orthogonal complement of $L^{(1)}$ in $L$ which is the sum of
irreducible sublattices of $L$ of {\it E}-type.
Let $L^{(2)}=L^{(2,1)}\oplus \cds \oplus L^{(2,n)}$ be the decomposition of
$L^{(2)}$ into irreducible components.
Then there exists $N\in \N$ such that we can embed $L^{(1)}$ into $D_{2N}$ and
$L^{(2,i)}$ into $E_8$ for $1\leq i\leq n$.
Set $R^0=D_{2N}$ and $R^i=E_8$ for $1\leq i\leq n$.
By using the embeddings above, we identify $L$ as a sublattice of
$R^0\oplus R^1\oplus \cds \oplus R^n \simeq D_{2N}\oplus E_8^{\oplus n}$
by which we regard $V_{\sqrt{2}L}^+$ as a sub VOA of $V_{\sqrt{2}(R^0\oplus \cds \oplus R^n)}^+$.
We shall show that if $e$ is an Ising vector of $V_{\sqrt{2}(R^0\oplus \cds \oplus R^n)}^+$,
then $e$ is contained in $V_{\sqrt{2}R^0}^+\tensor \cds \tensor V_{\sqrt{2}R^n}^+$.
Recall the code $\mathcal{C}_N$ of length $4N$ defined by \eqref{eq:c_n}.
Let $C$ be the dual code of the code generated by
$\mathcal{C}_N\oplus D_0(H_8)^{\oplus n}$ and $\gamma=(1010\dots 10)\in \Z_2^{4N+16n}$.
As shown in Sec.\ \ref{sec:4.1}, $V_{\sqrt{2}(R^0\oplus \cds \oplus R^n)}^+$ is
isomorphic to a code VOA $V_C$ associated to the code $C$.
Then by Proposition \ref{prop:3.7} there exists an embedding
$\phi:H_8 \hookrightarrow C$ such that $e\in V_{\phi(H_8)} \subset V_C$.
It is not difficult to verify that $\supp(\phi(H_8))$ is contained in one of the
direct summands of $\mathcal{C}_n^\perp \oplus (D_0(H_8)^\perp)^{\oplus n}$.
This implies that $e$ must be in $V_{\sqrt{2}R^0}^+\tensor \cds \tensor V_{\sqrt{2}R^n}^+$
as we claimed.
Then $e\in V_{\sqrt{2}R^i}^+$ for some $0\leq i\leq n$.
Therefore, by Theorem \ref{thm:4.4}, either there is a root
$\alpha \in R^0\oplus \cds \oplus R^n$ such that $e$ is of the form $w^\pm(\alpha)$ as in
\eqref{eq:4.3}, or $e\in V_{\sqrt{2}R^i}^+$ for some  $i>0$ and $e$ is of the form
\eqref{eq:4.5} in $V_{\sqrt{2}R^i}^+$.
Anyway, there is a sublattice $K$ of $R^0\oplus \cds \oplus R^n$ isometric to
$A_1$ or $E_8$ such that
$e\in V_{\sqrt{2}K}^+\subset V_{\sqrt{2}(R^0\oplus \cds \oplus R^n)}^+$.
If $e$ is taken from $V_{\sqrt{2}L}^+$, then $K\subset L$ as we know the explicit form
of $e$ in $V_{\sqrt{2}(R^0\oplus \cds \oplus R^n)}^+$.
This completes the proof.
\qed

\begin{rem}\label{rem:4.7}
  With reference to the theorem above,
  we believe in that the same is true for the $\Z_2$-orbifold $V_L^+$ of a lattice
  VOA associated with any even lattice $L$ without roots.
  For example, it is shown by Shimakura \cite{S2} that this is true for the
  $\Z_2$-orbifold $V_\Lambda^+$ of the Leech Lattice VOA $V_\Lambda$.
  However, we do not know the answer at present.
\end{rem}

\subsection{Ising vectors of $M_R$}

Let $R$ be an indecomposable root lattice as before. We shall
determine the Ising vectors of $M_R$. Recall that
\[
 w^\pm(\alpha)= \fr{1}{8}\alpha_{(-1)}^2\vac
    \pm \fr{1}{4}\l( e^{\sqrt{2}\alpha}+e^{-\sqrt{2}\alpha}\r),
    \quad \text{ for } \al\in R,
\]
and
\[
 \tilde{\w}=\tilde{\w}_R=\fr{2}{h+2}\w+\fr{1}{h+2}\sum_{\alpha \in \Phi(R)}
  e^{\sqrt{2} \alpha},
 \]
where $h$ denotes the Coxeter number of $\Phi(R)$ as before. By
definition, it is easy to show the following.

\begin{lem}\label{lem:4.7}
  Let $R$ be a root lattice and consider $V_{\sqrt{2}R}^+$.
  Let $\alpha,\beta\in \Phi(R)$.
  Then
  \\
  (1) $\la w^+(\alpha),w^-(\beta)\ra= 1/32$ if $\la \alpha,\beta \ra=\pm 1$ and
      $0$ otherwise.
  \\
  (2) $\la w^\pm (\alpha),w^\pm(\beta)\ra= 1/4$ if $\alpha =\pm \beta$, $1/32$
      if $\la \alpha,\beta\ra =\pm 1$ and $0$ if $\la \alpha,\beta\ra =0$.
  \\
  (3) $\la w^-(\alpha), \tilde{\w}\ra =0$ and $\la w^+(\alpha), \tilde{\w}\ra =1/(h+2)$,
      where $h$ is the Coxeter number of $\Phi(R)$.
  \\
  (4) For $R=E_8$, let $x\in E_8$.
      Then we have $\la \varphi_x\tilde{\w},\tilde{\w}\ra=1/4$ if $x\in 2E_8$,
      $1/32$ if $x$ is represented   by a root in $E_8/2E_8$ and $0$ otherwise.
\end{lem}

\begin{lem}\label{lem:4.8}
  Let $R$ be an indecomposable root lattice.
  The weight two subspace of $M_R$ is spanned by
  $\{ w^-(\alpha) \mid \alpha\in \Phi(R)\}$.
\end{lem}

\pf
We note that $M_R=\ker_{V_{\sqrt{2}R}^+}(\tilde{\w}_{(1)})$ since
$V_{\sqrt{2}R}^+$ is a completely reducible $\vir(\tilde{\w})$-module.
It is clear that $\{ w^-(\alpha) \mid \alpha \in \Phi(R)\}$ is a set of
linearly independent vectors of the weight two subspace of $M_R$, thanks to
(3) of Lemma \ref{lem:4.7}.
If $R=A_n$, it is shown in \cite{DLMN} that the set
$\{ w^\pm(\alpha) \mid \alpha \in \Phi(R)\}$ is a basis of the
weight two subspace of $V_{\sqrt{2}R}^+$.
Therefore, by (3) of Lemma \ref{lem:4.7}, the assertion holds in this case.
It is shown in \cite{DLY} that $M_{D_n}\simeq V_{\sqrt{2}A_{n-1}}^+$.
By Theorem \ref{thm:4.4} and (3) of Lemma \ref{lem:4.7}, the assertion also
holds if $R=D_n$.

Consider the case that $R=E_6$ or $E_7$.
The vacuum character of $M_R$ is computed in Sec.\ \ref{sec:6}
and one can verify that $\dim (M_R)_2= \abs{\Phi(R)}/2$.
So the assertion follows.
If $R=E_8$, then $M_{E_8}$ is a code VOA and it is easy to know the
structure code from which we can compute the vacuum character.
As a result, one also has $\dim (M_{E_8})_2=\abs{\Phi(E_8)}/2$.
This completes the proof.
\qed

\section{Automorphism group of $M_R$}

We will determine the automorphism group of the commutant subalgebra $M_R$ where
$R$ is an indecomposable root lattice.
It is shown that $\aut(M_R)$ always contains a half of the Weyl group $W(R)$ of $R$.

\subsection{$\mathbf{Aut}(M_R)$: the case $R\ne E_8$}

Let $R$ be an indecomposable root lattice.
We suppose that $R$ is not the $E_8$-lattice.
In this case, due to Theorem \ref{thm:4.4} and Lemma \ref{lem:4.7}, the set of
Ising vectors of $M_R$ is given by $E_R:=\{ w^-(\alpha)\mid \alpha \in \Phi(R)\}$.
We have shown in Lemma \ref{lem:4.8} that the weight two subspace of $M_R$ is
spanned by $E_R$.
In fact, $M_R$ is generated by $E_R$ as a VOA.

\begin{prop}\label{prop:key}
  $M_R$ is generated by its weight two subspace.
\end{prop}

\pf
The assertion is already shown in \cite{LS} if $R=A_n$.
Since $M_{D_n}\simeq V_{\sqrt{2}A_{n-1}}^+$ by \cite{DLY}, the assertion
is also true if $R=D_n$.
So we only need to show the cases for $R=E_6$ and $R=E_7$.
The proof is rather technical so that it will be given in Section \ref{sec:6}.
\qed
\vsb

Denote by $\Omega(E_R)$ the permutation group on $E_R$ and we define
$$
  \aut(E_R):=\{ \rho \in \Omega(E_R) \mid \la \rho e,\rho f\ra=\la e,f\ra
  \text{ for all } e,f\in E_R\} .
$$
Since $M_R$ is generated by $E_R$ as a VOA, by restriction map we have an injection
from $\aut(M_R)$ to $\aut(E_R)$.
On the other hand, $g\in \aut(R)$ acts on $E_R$ by $gw^-(\alpha):=w^-(g\alpha)$.
Hence we also have  a group homomorphism $\phi:\aut(R)\to \aut(E_R)$.

\begin{lem}\label{lem:5.2}
  $\phi$ is surjective.
\end{lem}

\pf
Let $\rho \in \aut(E_R)$.
Take a simple system $\Delta=\{ \alpha_1,\dots,\alpha_\ell\}$ of $\Phi(R)$, where
$\ell$ is the rank of $R$.
By Proposition \ref{prop:3.3} and Lemma \ref{lem:4.7}, one can easily verify that
$\{ w^-(\alpha_1),\dots,w^-(\alpha_\ell)\}$ is a set of generators of the Griess algebra
of $M_R$.
Write $\rho w^-(\alpha_i)=w^-(\beta_i)$ for $i=1,\dots,\ell$.
By Proposition \ref{prop:3.3} and Lemma \ref{lem:4.7}, there is a simple system $\Delta'$
of $\Phi(R)$ such that $\{ \pm \beta_1,\dots,\pm\beta_\ell\}=\Delta'\cup (-\Delta')$.
Therefore, by choosing suitable signs, we can take the representatives
$\{ \beta_i \mid 1\leq i\leq \ell\}$ to be a simple system of $\Phi(R)$.
Then we can find $g\in \aut(R)$ such that $g\alpha_i=\beta_i$ for $i=1,\dots,\ell$.
This implies that $\rho=g$ since $\{ w^-(\alpha_i) \mid 1\leq i\leq \ell\}$ is a set
of generators of $(M_R)_2=\Span_\C E_R$.
This completes the proof.
\qed
\vsb

By definition, it is clear that $\ker \phi= \la \pm 1\ra$.
Therefore, we have:

\begin{prop}\label{prop:5.3}
  $\aut (E_R)\simeq \aut(R)/\la \pm 1\ra$.
\end{prop}

The Weyl group $W(R)$ also acts on $E_R$.
Let us denote by $r_\alpha$ the reflection on $R$ defined by a root $\alpha\in \Phi(R)$.
It follows from Proposition \ref{prop:3.3} and Lemma \ref{lem:4.7} that
$\sigma_{w^-(\alpha)} w^-(\beta)=w^-(r_\alpha \beta)$ for $\alpha,\beta\in \Phi(R)$.
Therefore, the group generated by $\{ \sigma_e \mid e\in E_R\}$ realizes the
action of the Weyl group on $E_R$.
The shapes of the Weyl group $W(R)$ and the automorphism group $\aut(R)$ are
as follows (cf.\ \cite{CS}).

\renewcommand{\arraystretch}{1.5}
\begin{center}
\begin{tabular}{|c||c|c|}\hline
  $\q R\q $ & $\qq W(R)\qq$ & $\qq \aut(R)\qq$ \\ \hhline{|=#=|=|}
  $A_n$ & $\mathrm{S}_{n+1}$ & $\mathrm{S}_{n+1}\times 2$ \\\hline
  $D_4$ & $2^3:\mathrm{S}_4$ & $(2^3:\mathrm{S}_4):\mathrm{S}_3$ \\\hline
  $D_n$ ($n>4$) & $2^{n-1}:\mathrm{S}_n$ & $2^{n-1}:\mathrm{S}_n:2$ \\\hline
  $E_6$ & $\mathrm{U}_4(2):2$ & $2.\mathrm{U}_4(2):2$ \\\hline
  $E_7$ & $2\times \mathrm{Sp}_6(2)$ & $2\times \mathrm{Sp}_6(2)$\\\hline
\end{tabular}
\vsb
\end{center}

\begin{thm}\label{thm:5.4}
  Suppose $R\ne E_8$.
  Then $\aut(M_R)$ is as follows.
  \begin{center}
  \begin{tabular}{|c||c|c|c|c|c|c|}\hline
    $R$ & $A_1$ & $A_n$ $(n>1)$ & $D_4$ & $D_n$ $(n>4)$ & $E_6$ & $E_7$ \\ \hline
    $\aut(M_R)$
    & $1$
    & $\mathrm{S}_{n+1}$
    & $(2^2:\mathrm{S}_4):\mathrm{S}_3$
    & $2^{n-1}:\mathrm{S}_n$
    & $\mathrm{U}_4(2):2\simeq \mathrm{O}^-_6(2)$
    & $\mathrm{Sp}_6(2)$
    \\ \hline
  \end{tabular}
  \end{center}
\end{thm}

\pf
Let $R$ be either $A_n$, $n>1$, $E_6$ or $E_7$.
Under the injection, we can consider $\aut(M_R)$ as a subgroup of $\aut(E_R)$.
By Proposition \ref{prop:5.3}, we have the following relation.
$$
  \aut(M_R)\hookrightarrow \aut(E_R)\simeq \aut(R)/\la \pm 1\ra .
$$
Since the action of the Weyl group on $E_R$ is realized by $\sigma$-type involutions
$\{ \sigma_e \mid e\in E_R\}$, the subgroup of $\aut(M_R)$ generated by
$\{ \sigma_e \mid e\in E_R\}$ coincides with $\aut(R)/\la \pm 1\ra$.
Hence $\aut(M_R)\simeq \aut(R)/\la\pm 1\ra$.
This establishes the cases (i), (iii) and (iv).
The isomorphism in the case (ii) follows from the isomorphism
$M_{D_n}\simeq V_{\sqrt{2}A_{n-1}}^+$ shown in \cite{DLY} and the description of
$\aut(V_{\sqrt{2}A_n}^+)$ obtained in \cite{S1}.
The case $R=A_1$ is trivial since $M_{A_1}\simeq L(\shf,0)$.
This completes the proof.
\qed

\subsection{$\mathbf{Aut} (M_{E_8})\simeq \mathrm{Sp}_8(2)$}

By Theorem \ref{thm:4.4}, the set of Ising vectors of $V_{\sqrt{2}E_8}^+$ is given by
\begin{equation}\label{eq:4.6}
  \{ w^\pm(\alpha) \mid \alpha \in \Phi(E_8)\}
  \cup \{ \varphi_x \tilde{\w} \mid x\in E_8\} .
\end{equation}

\begin{lem}\label{lem:5.5}
  Let $R$ be an indecomposable root lattice.
  Then all the Ising vectors of $V_{\sqrt{2}R}^+$ are conjugate under
  $\aut (V_{\sqrt{2}R}^+)$.
\end{lem}

\pf
If $R\ne E_8$, then the set of Ising vectors of $V_{\sqrt{2}R}^+$ is provided by
$\{ w^\pm(\alpha) \mid \alpha \in \Phi(R)\}$.
It is clear that $V_{\sqrt{2}R}^+$ affords an action of the Weyl group $W(R)$
associated to the root system $\Phi(R)$ and $W(R)$ transitively acts on both
$\{ w^+(\alpha) \mid \alpha\in \Phi(R)\}$ and $\{ w^-(\beta)\mid \beta\in \Phi(R)\}$
as $R$ is indecomposable.
By (2) of Proposition \ref{prop:3.3} and (1) of Lemma \ref{lem:4.7}, there is a pair
$\{ w^+(\alpha),w^-(\beta)\}$ of Ising vectors which are conjugate under
$\aut (V_{\sqrt{2}R}^+)$.
Therefore, all the Ising vectors are conjugate.

Now assume that $R=E_8$.
We have seen that $\aut (V_{\sqrt{2}E_8}^+)$ acts on
$\{ w^\pm(\alpha) \mid \alpha \in \Phi(E_8)\}$ transitively.
It is obvious that Ising vectors of the form \eqref{eq:4.5} are mutually conjugate
under $\aut (V_{\sqrt{2}E_8}^+)$.
Then again by (2) of Proposition \ref{prop:3.3} and (3) of Lemma \ref{lem:4.7},
all the Ising vectors of $V_{\sqrt{2}E_8}^+$ are conjugate.
\qed
\vsb

As an immediate consequence, we have:

\begin{cor}\label{cor:5.6}
  For any Ising vector $e$ of $V_{\sqrt{2}E_8}^+$,
  $\com_{V_{\sqrt{2}E_8}^+}(\vir(e))\simeq M_{E_8}$.
\end{cor}

\begin{lem}\label{lem:5.7}
  $M_{E_8}$ is generated by its weight two subspace.
\end{lem}

\pf
It is shown in Sec.\ \ref{sec:4.1} that $V_{\sqrt{2}E_8}^+$ is isomorphic to
the code VOA $V_{\mathrm{RM}(2,4)}$ associated to the Reed-Muller code $\mathrm{RM}(2,4)$.
By Corollary \ref{cor:5.6}, we may assume that $\tilde{\w}$ is the 1st summand of
the standard Ising frame of $V_{\mathrm{RM}(2,4)}$.
Then $M_{E_8}$ is also a code VOA and the associated code is obtained by collecting
all the vectors of $\mathrm{RM}(2,4)$ whose 1st entry is 0 and dropping their 1st entry 0.
It is easy to see that this code is generated by weight four vectors.
Therefore, $M_{E_8}$ is generated by its weight two subspace.
\qed

\begin{lem}\label{lem:5.8}
  (1) An Ising vector of $M_{E_8}$ is equal to either $w^-(\alpha)$, $\alpha \in \Phi(E_8)$
  or $\varphi_x\tilde{\w}$, $x\in E_8$ and $x+2E_8 \in E_8/2E_8$ is represented by a norm
  four vector of $E_8$.
  Hence, there are $\mathit{120+135=255}$ Ising vectors of $M_{E_8}$.
  \\
  (2) These 255 Ising vectors are mutually conjugate under $\aut(M_{E_8})$.
\end{lem}

\pf
(1) Since $V_{\sqrt{2}E_8}^+$ has a trivial weight one subspace, so does $M_{E_8}$.
By (1) of Proposition \ref{prop:3.3}, an Ising vector $e$ of $V_{\sqrt{2}E_8}^+$
belongs to the commutant subalgebra $M_{E_8}=\com_{V_{\sqrt{2}E_8}^+}(\vir(\tilde{\w}))$
if and only if $\la e,\tilde{\w}\ra =0$.
So the assertion follows from Lemma \ref{lem:4.7}.

(2) It is clear that $w^-(\alpha)$, $\alpha\in \Phi(E_8)$, are mutually conjugate
under $\aut(M_{E_8})$ since $M_{E_8}$ affords a natural action of the Weyl group $W(E_8)$.
One can also verify that for each $\varphi_x \tilde{\w}\in M_{E_8}$, there is an Ising
vector $w^-(\beta)\in M_{E_8}$ such that $\la \varphi_x\tilde{\w},w^-(\beta)\ra =1/32$.
Therefore, all the Ising vectors are mutually conjugate by (2) of Proposition \ref{prop:3.3}.
\qed

\begin{thm}\label{thm:5.9}
  $\aut(M_{E_8})\simeq \mathrm{Sp}_8(2)$.
\end{thm}

\pf
Denote by $E$ the set of Ising vectors of $M_{E_8}$ and set
$$
  \sigma(E):=\{\sigma_e \in \aut(M_{E_8}) \mid e\in E\} .
$$
Let $G$ be the 3-transposition subgroup of $\aut(M_{E_8})$
generated by $\sigma(E)$. It is shown in \cite{KM,Ma} that $G$ is
isomorphic to a simple group $\mathrm{Sp}_8(2)$. Since $E$ is
invariant under $\aut(M_{E_8})$, $\sigma(E)$ is a normal set of
$\aut(M_{E_8})$ and hence $G$ is a normal subgroup of
$\aut(M_{E_8})$. It is shown in \cite{ATLAS} that
$\aut(G)=\mathrm{Inn}(G)\simeq G$. So the kernel of the conjugate
action of $\aut(M_{E_8})$ on $G$ is $C_{\aut(M_{E_8})}(G)$ and
$\aut(M_{E_8})\simeq G \times C_{\aut(M_{E_8})}(G)$. Let $g\in
C_{\aut(M_{E_8})}(G)$. Since $g$ commutes with $\sigma(E)$, $g$
acts on $E$ identically by Lemma \ref{lem:3.13}. By Lemma
\ref{lem:4.8}, $g$ is also identical on the weight two subspace of
$M_{E_8}$. Since $M_{E_8}$ is generated by $E$ as a VOA by Lemma
\ref{lem:5.7}, $g$ is trivial on $M_{E_8}$. Thus
$C_{\aut(M_{E_8})}(G)=1$ and $\aut(M_{E_8})=G\simeq \mathrm{Sp}_8(2)$.
\qed


\section{3-transposition group and inductive structure}

We will study a relation between an inductive structure of a 3-transposition group
acting on a vertex operator algebra and its commutant subalgebra.
As an example, we study a commutant subalgebra of $V_{\sqrt{2}E_8}^+$ having
$\mathrm{O}_8^-(2)$ as its full automorphism group, which corresponds
to an inductive structure $\mathrm{O}_{10}^+(2)^{(2)}\simeq \mathrm{O}_8^-(2)$.
Note that $\aut(V_{\sqrt{2}E_8}^+)\simeq \mathrm{O}_{10}^+(2)$  by \cite{G1,S1}.

\subsection{3-transposition group}

\begin{df}\label{df:5.1}
  A (finite) {\it 3-transposition group} is a pair $(G,\mathcal{D})$ of a finite group
  $G$ and a normal set $\mathcal{D}$ of involutions in $G$ such that $G$ is generated by
  $\mathcal{D}$ and if $x,y\in \mathcal{D}$ then the order of $xy$ is either 1,2 or 3.
\end{df}

A {\it partial linear space} consists of a set $X$ called the set of points and a set
$\mathcal{L}$ of subsets of $X$ called the set of lines such that any two
points lie on  at most one line and any line has at least two points.

To a 3-transposition group $(G,\mathcal{D})$, we can associate a partial linear
space called the {\it Fischer space} $(X,\mathcal{L})$ of $(G,\mathcal{D})$
as follows (cf.\ \cite{Ma} and references therein).
The set of points $X$ is $\mathcal{D}$ and the set of lines $\mathcal{L}$ is
such that a subset $\ell\subset \mathcal{D}$ is a line if and only if $\ell$
consists of three points which generate a subgroup of $G$ isomorphic to $\mathrm{S}_3$.
For $x,y\in \mathcal{D}$, we write $x\sim y$ if they lie on a line and
$x\perp y$ if not.

The {\it collinearity graph}  of $(X,\mathcal{L})$ is a graph $\Gamma$
whose vertex set is $X$ and $x,y\in X$ are adjacent if $x\sim y$, i.e.,
they are incident to a common line.

A 3-transposition group $(G,D)$ is referred to as {\it of symplectic type} if
the affine plane of order 3 does not occur in the associated Fischer space.
It is called {\it indecomposable} if the associated Fischer space has the connected
collinearity graph.

\begin{rem}
  Let $\ell_1$ and $\ell_2$ be two distinct lines with a common point in a
  Fischer space.
  It is known that the subspace generated by $\ell_1$ and $\ell_2$ is either
  isomorphic to a dual affine plane of order $2$ or an affine plane of order 3
  (cf.\ \cite{Asc}).
  Thus if $(G,D)$ is of symplectic type, only the dual affine plane of order
  $2$ can occur in the associated Fischer space.
  In this case, the subgroup generated by the corresponding involutions in $\ell_1$
  and $\ell_2$ will be isomorphic to the symmetry group $\mathrm{S}_4$.
\end{rem}

Let $(G,\mathcal{D})$ be a 3-transposition group and $(X,\mathcal{L})$ the
associated Fischer space.
Define $\sigma: X\to \aut (X)$ by $\sigma_x(y):=y^x$ for $x\in X$, and for
$x_1,x_2,\dots \in X$, denote by
$$
  X_{x_1,x_2,\dots}:=\{ x\in X\mid x\perp x_i,\ i=1,2,\dots\} .
$$
If $G$ is indecomposable, $\sigma(G)$ acts on $X=\mathcal{D}$ transitively.
For, if $x\sim y$ with $x,y\in X$, then there exists $z\in X$ such that
$x\sim z\sim y$, which implies $\sigma_x(y)=\sigma_y(x)=z$.
Thus $x$ and $y$ are conjugate under the action of $\sigma(G)$.
Repeating this, we can establish the transitivity.

Given an (indecomposable) 3-transposition group $(G,D)$, one can consider
its {\it inductive structure} as follows.
Let $D^{(1)}$ be the set of elements of $D$ which commute with a fixed element
of $D$ and $D^{(2)}$ be that with two fixed non-commuting elements of $D$.
We set $G^{(1)}=\la D^{(1)}\ra$ and $G^{(2)}=\la D^{(2)}\ra$.
The group structure $(G^{(i)},D^{(i)})$ is called the inductive structure of
$(G,D)$.

\begin{prop}\label{prop:7.3}
  Let $(G,\mathcal{D})$ be a centerfree indecomposable 3-transposition group
  of symplectic type and $(X,\mathcal{L})$ the associated Fischer space.
  Suppose that $X_{x,y}\ne \emptyset$ for any distinct $x,y\in X$.
  If $C_{\aut(X)}(\sigma_x)|_{X_x}=C_G(\sigma_x)|_{X_x}$ for any $x\in X$,
  then $\aut(X)=G$.
\end{prop}

\pf
Let $g\in \aut(X)$. Take any $x\in X=\mathcal{D}$.
By the transitivity, there exists $\rho_1\in G$ such that $g(x)=\rho_1(x)$.
Hence $\rho_1^{-1}g|_{X_x}\in C_{\aut(X)}(\sigma_x)|_{X_x}=C_G(\sigma_x)|_{X_x}$.
Then there exists $\rho_2\in G$ such that $\rho_2^{-1}\rho_1^{-1}g=\id$ on
$X_x\cup \{ x\}$.
Take $y\in X$ such that $y\sim x$.
Set $\rho=\rho_1\rho_2$.
We show that $\{ x,y,y^x\}=X_{x,y}^\perp=\{ z \in X \mid z \perp X_{x,y}\}$,
from which we deduce that $\rho^{-1}g$ stabilizes the set $\{ x,y,y^x\}$.
Assume that there exists $z\in X_{x,y}^\perp\setminus \{ x,y,x^y\}$.
We may assume that $z\sim x$.
Since $(X,\mathcal{L})$ is a Fischer space of symplectic type,
the subspace $\la x,y,z\ra$ generated by $x,y$ and $z$ is isomorphic to
a dual affine plane of order $2$ and the configuration in the Fischer space
is as in Figure 1:
\vsb\\
\unitlength 0.1in
\begin{picture}( 39.8500, 18.5000)(  4.0000,-18.7000)
%
\special{pn 20}%
\special{sh 1}%
\special{ar 2332 1670 10 10 0  6.28318530717959E+0000}%
\special{sh 1}%
\special{ar 3310 190 10 10 0  6.28318530717959E+0000}%
\special{sh 1}%
\special{ar 4288 1670 10 10 0  6.28318530717959E+0000}%
\special{sh 1}%
\special{ar 3868 1036 10 10 0  6.28318530717959E+0000}%
\special{sh 1}%
\special{ar 2752 1036 10 10 0  6.28318530717959E+0000}%
\special{sh 1}%
\special{ar 2752 1036 10 10 0  6.28318530717959E+0000}%
%
\special{pn 8}%
\special{pa 2332 1670}%
\special{pa 3310 190}%
\special{fp}%
\special{pa 3310 190}%
\special{pa 4288 1670}%
\special{fp}%
\special{pa 4288 1670}%
\special{pa 2752 1036}%
\special{fp}%
\special{pa 2332 1670}%
\special{pa 3868 1036}%
\special{fp}%
\put(33.7900,-1.9000){\makebox(0,0)[lb]{$x$}}%
\put(39.6600,-10.3600){\makebox(0,0)[lb]{$z$}}%
\put(43.8500,-16.7000){\makebox(0,0)[lb]{$z^x=x^z$ }}%
\put(26.4600,-10.3600){\makebox(0,0)[rb]{$y$}}%
\put(22.9000,-16.7000){\makebox(0,0)[rb]{$y^x=x^y$}}%
%
\special{pn 20}%
\special{sh 1}%
\special{ar 3310 1268 10 10 0  6.28318530717959E+0000}%
\special{sh 1}%
\special{ar 3310 1268 10 10 0  6.28318530717959E+0000}%
\put(30.8000,-20.4000){\makebox(0,0)[lb]{Fig. 1}}%
\put(31.4000,-15.6000){\makebox(0,0)[lb]{$(y^x)^z$}}%
\put(4.0000,-8.7000){\makebox(0,0)[lb]{\quad\quad}}%
\end{picture}%
\vsv\\
Since $X$ is indecomposable  and $X_{x,y}\ne \emptyset$, there exists
$w\in X_{x,y}$ such that $w$ is collinear to a point $p$ in $\la x,y,z\ra$.
By Fig.\ 1., $p=z$, $x^z$ or $(y^x)^z$.
Nevertheless, $z\in X_{x,y}^\perp$.
Hence we have $p=x^z$ or $p=(y^x)^z$.
In either case, $p\sim z$ and $z^p=x$ or $y^x$.
Now consider the subspace generated by $w,z$ and $p$.
Then we have the configuration as in Figure 2:
\vsb\\
\unitlength 0.1in
\begin{picture}( 39.8500, 18.5000)(  4.0000,-18.7000)
%
\special{pn 20}%
\special{sh 1}%
\special{ar 2332 1670 10 10 0  6.28318530717959E+0000}%
\special{sh 1}%
\special{ar 3310 190 10 10 0  6.28318530717959E+0000}%
\special{sh 1}%
\special{ar 4288 1670 10 10 0  6.28318530717959E+0000}%
\special{sh 1}%
\special{ar 3868 1036 10 10 0  6.28318530717959E+0000}%
\special{sh 1}%
\special{ar 2752 1036 10 10 0  6.28318530717959E+0000}%
\special{sh 1}%
\special{ar 2752 1036 10 10 0  6.28318530717959E+0000}%
%
\special{pn 8}%
\special{pa 2332 1670}%
\special{pa 3310 190}%
\special{fp}%
\special{pa 3310 190}%
\special{pa 4288 1670}%
\special{fp}%
\special{pa 4288 1670}%
\special{pa 2752 1036}%
\special{fp}%
\special{pa 2332 1670}%
\special{pa 3868 1036}%
\special{fp}%
\put(33.7900,-1.9000){\makebox(0,0)[lb]{$p$}}%
\put(39.6600,-10.3600){\makebox(0,0)[lb]{$z$}}%
\put(43.8500,-16.7000){\makebox(0,0)[lb]{$z^p=p^z$ }}%
\put(26.4600,-10.3600){\makebox(0,0)[rb]{$w$}}%
\put(22.9000,-16.7000){\makebox(0,0)[rb]{$w^p=p^w$}}%
%
\special{pn 20}%
\special{sh 1}%
\special{ar 3310 1268 10 10 0  6.28318530717959E+0000}%
\special{sh 1}%
\special{ar 3310 1268 10 10 0  6.28318530717959E+0000}%
\put(30.8000,-20.4000){\makebox(0,0)[lb]{Fig. 2}}%
\put(31.4000,-15.6000){\makebox(0,0)[lb]{$(w^p)^z$}}%
\put(4.0000,-8.7000){\makebox(0,0)[lb]{\quad\quad}}%
\end{picture}%
\vsv\\
In this case, $w\sim z^p \in \{x,y^x\}$ which contradicts the fact that
$w\in X_{x,y}$.
Hence $X_{x,y}^\perp=\{ x,y,x^y\}$.
Thus by replacing $\rho$ by $\rho \sigma_x$ if needed, we may assume that
$\rho^{-1}g(y)=y$. Then $\rho^{-1}g$ acts as identity on $X_x\cup \{ x,y,y^x\}$.
Now we claim that $\rho^{-1}g=\id$ on $X$.
Take any $t\in X$ such that $t\not\in X_x\cup \{x,y,y^x\}$.
Then $\{x,y,x^y\}$ and $\{x,t, x^t\}$ are two distinct lines.
Since $(X,\mathcal{L})$ is of symplectic type, the Fischer subspace generated by
$\{x,y, t\}$ has the configuration as in Fig.\ 3:
\vsb\\
\unitlength 0.1in
\begin{picture}( 39.8500, 18.5000)(  4.0000,-18.7000)
%
\special{pn 20}%
\special{sh 1}%
\special{ar 2332 1670 10 10 0  6.28318530717959E+0000}%
\special{sh 1}%
\special{ar 3310 190 10 10 0  6.28318530717959E+0000}%
\special{sh 1}%
\special{ar 4288 1670 10 10 0  6.28318530717959E+0000}%
\special{sh 1}%
\special{ar 3868 1036 10 10 0  6.28318530717959E+0000}%
\special{sh 1}%
\special{ar 2752 1036 10 10 0  6.28318530717959E+0000}%
\special{sh 1}%
\special{ar 2752 1036 10 10 0  6.28318530717959E+0000}%
%
\special{pn 8}%
\special{pa 2332 1670}%
\special{pa 3310 190}%
\special{fp}%
\special{pa 3310 190}%
\special{pa 4288 1670}%
\special{fp}%
\special{pa 4288 1670}%
\special{pa 2752 1036}%
\special{fp}%
\special{pa 2332 1670}%
\special{pa 3868 1036}%
\special{fp}%
\put(33.7900,-1.9000){\makebox(0,0)[lb]{$x$}}%
\put(39.6600,-10.3600){\makebox(0,0)[lb]{$t$}}%
\put(43.8500,-16.7000){\makebox(0,0)[lb]{$x^t=t^x$ }}%
\put(26.4600,-10.3600){\makebox(0,0)[rb]{$y$}}%
\put(22.9000,-16.7000){\makebox(0,0)[rb]{$y^x=x^y$}}%
%
\special{pn 20}%
\special{sh 1}%
\special{ar 3310 1268 10 10 0  6.28318530717959E+0000}%
\special{sh 1}%
\special{ar 3310 1268 10 10 0  6.28318530717959E+0000}%
\put(30.8000,-20.4000){\makebox(0,0)[lb]{Fig. 3}}%
\put(32.7000,-14.6000){\makebox(0,0)[lb]{$s$}}%
\put(4.0000,-8.7000){\makebox(0,0)[lb]{\quad\quad}}%
\end{picture}%
\vsv\\
Hence, there exists  $s\in X_x \cap \la x,y,t\ra$ such that $\la
x,y,t\ra = \la x,y,s\ra\simeq \mathrm{S}_4$ and $t\in \la x,y,s\ra$.
Since $\rho^{-1}g$ acts trivially on $\{ x,y,s\}$, it acts trivially on
$t$, also. Thus $G=\aut(X)$.
\qed

\subsection{Automorphism group of a commutant subalgebra}

We consider a correspondence between an inductive structure of a 3-transposition group
and a commutant subalgebra structure of a vertex operator algebra on which the group acts.
We study an inductive structure $\mathrm{O}^+_{10}(2)^{(2)}\simeq
\mathrm{O}^-_8(2)$.
It is shown in \cite{G1,S1} that $\aut(V_{\sqrt{2}E_8}^+)\simeq \mathrm{O}_{10}^+(2)$.
By this inductive structure, we can find a commutant subalgebra of
$V_{\sqrt{2}E_8}^+$ on which $\mathrm{O}^-_8(2)$ acts.
We will show that in this case the inductive structure completely determines
the full automorphism group of the commutant subalgebra.

Take a root $\alpha_0\in \Phi(E_8)$ and set $L_1:=\Z \alpha_0$ and
$L_2:= (\Z \alpha_0)^\perp =\{ \beta \in E_8 \mid \la \alpha_0,\beta\ra= 0\}$.
Then $L_1$ and $L_2$ are sublattices of $E_8$ isometric to $A_1$ and $E_7$, respectively.
Let $U$ be the subalgebra of
$V_{\sqrt{2}E_8}^+$ generated by $\tilde{\w}$ and $\varphi_{\alpha_0} \tilde{\w}$.
It is shown in \cite{M1,LY1} that $U$ is isomorphic to a simple VOA of the form
$$
  L(\hf,0)\tensor L(\fr{7}{10},0)\oplus L(\hf,\hf)\tensor L(\fr{7}{10},\fr{3}{2}).
$$
By a direct computation, one has
$$
  \tilde{\w}_{(1)}\varphi_{\alpha_0} \tilde{\w}
  =\dfr{1}{4}\l( \tilde{\w} +\varphi_{\alpha_0} \tilde{\w}-w^+(\alpha_0)\r)
$$
so that $\sigma_{\tilde{\w}}\varphi_{\alpha_0} \tilde{\w}
=\sigma_{\varphi_{\alpha_0} \tilde{\w}} \tilde{\w}=w^+(\alpha_0)$
(cf.\ Proposition \ref{prop:3.3}).
Therefore, $U$ contains exactly three Ising vectors $\tilde{\w}$,
$\varphi_{\alpha_0} \tilde{\w}$ and $w^+(\alpha_0)$.
Set
\[
  U^c=\com_{V_{\sqrt{2}E_8}}(U).
\]
We shall show that $U^c$ is again generated by a set of Ising vectors and
its full automorphism group $\aut(U^c)$ is isomorphic to $\mathrm{O}_8^-(2)$.

\begin{prop}\label{prop:6.5+}
  There are exactly $136$ Ising vectors in the commutant subalgebra
  $U^c$ and they are mutually conjugate under the action of $\aut(U^c)$.
\end{prop}

\proof
Since $U$ is generated by $\tilde{\w}$ and $\varphi_{\alpha_0} \w$,
we have
$$
  U^c=\com_{V_{\sqrt{2}E_8}^+}\! \l(\vir(\w)\r)\cap
      \com_{V_{\sqrt{2}E_8}^+}\! \l(\vir(\varphi_{\alpha_0} \tilde{\w})\r) .
$$
By Proposition \ref{prop:4.3} and Lemma \ref{lem:4.7}, if $e$ is
an Ising vector in $U^c$, then either
\vsb

(1) $e= w^-(\be)$ for some $\be $ with $\la \al_0, \be \ra=0,\pm 2$\quad or

(2) $e=\varphi_x( \tilde{\om})$, $x+ 2E_8$ is represented by a norm 4 vector
    and $\la \al_0, x \ra \equiv 1 \mod 2$.
\vsb

\noindent
Note that the lattice $L_2=(\Z \alpha_0)^\perp$ is isomorphic to $E_7$, which has
$63$ positive roots.
Moreover, $\la \al_0, \be\ra =\pm 2$  if and only if $\be=\pm \al_0$.
Hence, there are exactly $1+63=64$ Ising vectors of the form $w^-(\be)$ in
the case (1).

For the case (2), there are exactly $128$ cosets $x+2E_8$ of $E_8/2E_8$
such that $\la \al_0, x \ra =1 \mod 2$.
Among them, $120 -64=56$ classes are represented by roots.
Hence there are $72$ Ising vectors of the form $\varphi_x (\tilde{\om})$ in $U^c$.
Therefore, there are totally $64+72=136$ Ising vectors in $U^c$.

Next we shall show all Ising vectors are mutually conjugate in $U^c$.
First we shall note that the Weyl group $W(L_2)\simeq W(E_7)$ is naturally a subgroup
of $\aut (U^c)$ since both $\tilde{\om}$ and $w^+(\al_0)$ are fixed by $W(L_2)$.
Hence, it is clear that $w^-(\be)$ and $w^-(\gamma)$ are conjugate if
$\be,\gamma\in L_2$.
Let $e=\varphi_x(\tilde{\om})$ be an Ising vector of $U^c$.
Then $x+ 2E_8$ is represented a norm 4 vector and $\la \al_0, x \ra \equiv 1 \mod 2$.
Since $E_8$ is generated by its roots, we may set $x=\be_1+\be_2$ such that
$\be_1$ and $\be_2$ are roots in $E_8$ and $\la \be_1, \be_2\ra =0$,
$\la \al_0,\be_1 \ra =1$ and $\la \al_0, \be_2 \ra =0$.
In this case, $\be_2\in (\Z \alpha_0 +\Z \beta_1)^\perp\simeq E_6$ and hence
there exists a root $\gamma\in (\Z \alpha_0+\Z \beta_1)^\perp$ such that
$\la \be_1,\gamma \ra =\la \al_0,\gamma \ra =0$ and $\la \be_2, \gamma \ra =1$.
Therefore, we have $\la \varphi_x(\tilde{\om}), w^-(\gamma)\ra =1/32$ and thus
$\varphi_x(\tilde{\om})$ is conjugate to $w^-(\gamma)$ by (2) of Proposition
\ref{prop:3.3}.

Finally, for any $e=\varphi_x( \tilde{\om})\in U^c$, we have
$\la \al_0, x \ra \equiv 1 \mod 2$.
Hence,  we have $\la \varphi_x(\tilde{\om}), w^-(\al_0)\ra =1/32$ and
$\varphi_x(\tilde{\om})$ is conjugate to $w^-(\al_0)$ also.
\qed

\begin{lem}
  $U^c$ is generated by its Ising vectors.
\end{lem}

\proof
Denote
\[
  L=\sqrt{2} (L_1\oplus L_2)\simeq \sqrt{2}A_1 \oplus \sqrt{2}E_7.
\]
Since $|\sqrt{2}E_8/L|=2$, there is $\gamma\in \sqrt{2}E_8$ such
that $\sqrt{2}E_8= L \cup (\gamma+L) $ and hence
\[
  V_{\sqrt{2}E_8} =V_L \oplus V_{\gamma+L}\q \mbox{and} \q
  V_{\sqrt{2}E_8}^+=V_L^+\oplus V_{\gamma+L}^+.
\]
Note that the quotient group structure $\sqrt{2}E_8/L$ induces an
automorphism $\rho\in \aut(V_{\sqrt{2}E_8}^+)$ such that
$\rho|_{V_L^+}=1$ and $\rho|_{V_{\gamma+L}^+}=-1$.

By definition (cf.\ \eqref{eq:2.3}), it is easy to show
\[
  \tilde{\w}_{L_2}(=\tilde{\w}_{E_7})
  = \dfr{4}{5}(\tilde{\w}+\varphi_{\alpha_0}\tilde{\w})-\dfr{1}{5}w^+(\alpha_0).
\]
Therefore, $\tilde{\w}_{L_2}\in U$ and the Virasoro element of $U$
is an orthogonal sum of $w^+(\al_0)$ and $\tilde{\w}_{L_2}$.
Hence the commutant subalgebra $U^c$ can be defined as follows.
\[
  U^c= \{ v\in V_{\sqrt{2}E_8}^+\mid (\tilde{\om}_{L_2})_{(1)}
          v=w^+(\al_0)_{(1)}v=0\} .
\]
Let
\[
\begin{array}{lllll}
  M^0 &=& U^c\cap V_L^+
      &=& \{ v\in V_{L}^+ \mid (\tilde{\om}_{L_2})_{(1)} v=w^+(\al_0)_{(1)}v=0\} ,
  \\
  M^1 &=& U^c\cap V_{\gamma+L}^+
      &=& \{ v\in V_{\gamma+ L}^+ \mid (\tilde{\om}_{L_2})_{(1)} v=w^+(\al_0)_{(1)} v=0\} .
\end{array}
\]
Then we have $U^c= M^0\oplus M^1$. Moreover, the automorphism
$\rho$ induces a natural action on $U^c$ such that $\rho|_{M^0}=1$ and $\rho|_{M^1}=-1$.
Note that
\[
  M^0\simeq  L(1/2,0)\otimes M_{E_7}
\]
and  $M^1$ is an irreducible $M^0$-module in this case (cf.\ \cite{DM}).

Now by Lemma \ref{lem:4.7} and Proposition \ref{prop:key},
we know that $M_{E_7}$ is generated by the Ising vectors of the form
$ w^-(\be)$, $\be \in L_2\simeq E_7$ .
Hence, $M^0\simeq L(\shf,0)\otimes M_{E_7}$ is generated by $w^-(\al_0)$ and
$\{ w^-(\be) \mid  \be \in L_2\simeq E_7\}$.

Let $W$ be the sub VOA generated by the set of all Ising vectors of $U^c$.
Then  $M^0\subset W$.
Since $W$ also contains Ising vectors of the form $\varphi_x(\tilde{\om})$
which is not contained in $M^0$ (cf.\ Lemma  \ref{lem:4.7}),
we know that $W\neq M^0$.
Hence, $W=M^0\oplus M^1=U^c$ as desired.
\qed

\begin{thm}
  $\aut (U^c)\simeq \mathrm{O}^-_8(2)$.
\end{thm}

\pf
Let $E$ be the set of Ising vectors of $U^c$ and $G$ the 3-transposition
subgroup of $\aut(U^c)$ generated by involutions $\{ \sigma_e\mid e\in E\}$.
It is shown in \cite{Ma} that $G\simeq \mathrm{O}^-_8(2)$.
By the proof of the previous lemma, we have
$$
  \tilde{\w}_{L_2}
  = \dfr{4}{5}(\tilde{\w} + \varphi_{\alpha_0} \tilde{\w})
    -\dfr{1}{5}w^+(\alpha_0) \in U
$$
and $\com_{U^c}(\vir(w^-(\alpha_0)))\simeq M_{E_7}$.
By Theorem \ref{thm:5.4}, $\aut(M_{E_7})$ is isomorphic to $\mathrm{Sp}_6(2)$
which is generated by $\sigma$-type involutions associated to Ising vectors
of $M_{E_7}$.
Since the 3-transposition subgroup of $G$ generated by
$\sigma(E)_{\sigma_{w^-(\alpha)}} =\{ \sigma_{w^-(\beta)} \mid \beta \in L_2\}$
is isomorphic to $\mathrm{Sp}_6(2)$, we can apply Proposition \ref{prop:7.3} to $G$
and we conclude that $\aut (U^c)=G\simeq \mathrm{O}^-_8(2)$ as $U^c$ is
generated by $E$.
\qed


\section{Decomposition of $M_R$}\label{sec:6}

Next, we will complete the proof of Proposition \ref{prop:key}.
Since it requires the notion of $W$-algebras,  we will first review
some basic facts about $W$-algebras.

\subsection{Modules over $W$-algebra}

Let $\Lambda_0$ and $\Lambda_1$ be the fundamental weights of the affine Lie algebra
$\hat{\mathfrak{sl}}_2(\C)$.
For positive integers $\ell, j$ with $0\leq j\leq \ell$, consider the irreducible
highest weight module $\eL (\ell,j)$ over $\hat{\mathfrak{sl}}_2(\C)$ with highest
weight $(\ell -j)\Lambda_0+j\Lambda_1$.
It is well-known that $\eL (\ell,0)$ forms a simple VOA and the integrable
$\hat{\mathfrak{sl}}_2(\C)$-modules $\eL (\ell,j)$, $0\leq j\leq \ell$,
provide all the inequivalent irreducible $\eL (\ell,0)$-modules (cf.\ \cite{FZ}).

Now let
$A_1^{\oplus \ell} = \Z \epsilon_1 \oplus \cdots \oplus \Z \epsilon_\ell$ be an
even lattice with $\la \epsilon_i,\epsilon_j\ra =2\delta_{i,j}$ and
$V_{A_1^{\oplus \ell}}$ the lattice VOA associated with $A_1^{\oplus\ell}$.
Then $V_{A_1^{\oplus \ell}}\simeq (V_{A_1})^{\tensor \ell} \simeq \eL (1,0)^{\tensor \ell}$.
Recall that the weight one subspace of $V_{A_1^{\oplus \ell}}$ forms a Lie
algebra by the Lie bracket $[x,y]:=x_{(0)} y$ for $x,y \in (V_{A_1^{\oplus \ell}})_1$.
Set
$H^{(\ell)}=({\epsilon_1} +\cdots + \epsilon_\ell)_{(-1)}\vac$,
$E^{(\ell)}=e^{\epsilon_1}+\cdots +e^{\epsilon_\ell}$ and
$F^{(\ell)}=e^{-\epsilon_1}+\cdots +e^{-\epsilon_\ell}$.
Then the subspace $\C H^{(\ell)}+\C E^{(\ell)}+\C F^{(\ell)}$ of the weight one
subspace of $V_{A_1^{\oplus \ell}}$ forms a simple Lie subalgebra isomorphic to
$\mathfrak{sl}_2(\C)$ and the sub VOA generated by
$\{ H^{(\ell)}, E^{(\ell)}, F^{(\ell)} \}$ is isomorphic to level $\ell$ affine
VOA $\eL (\ell,0)$ (cf.\ \cite{DL}).

Let $\gamma= \epsilon_1+\cdots+\epsilon_\ell\in A_1^{\oplus \ell}$.
Then $\gamma_{(-1)}\vac = H^{(\ell)}$ and it is easy to verify that
$$
  e^{\gamma}=\frac{1}{\ell!}(E^{(\ell)}_{(-1)})^{\ell-1}E^{(\ell)}.
$$
Thus $\eL (\ell,0)$ contains a subalgebra isomorphic to the lattice
VOA $V_{\Z \gamma}$.
We consider the following commutant subalgebra:
\begin{equation}\label{eq:6.1}
  W_\ell:= \com_{\eL (\ell,0)}(V_{\Z \gamma}).
\end{equation}
Note that the central of $W_\ell$ is equal to $3\ell/(\ell+2)-1=2(\ell-1)/(\ell+2)$.
It is clear that $\eL (\ell,0)$ contains a full sub VOA isomorphic to
$V_{\Z \gamma}\tensor W_\ell$ so that we can regard $\eL (\ell,j)$ as a
$V_{\Z\gamma}\tensor W_\ell$-module.
Following \cite{DL,Li1}, we introduce the following spaces:
\begin{equation}\label{eq:6.2}
  W_\ell(j,k)
  := \hom_{V_{\Z \gamma}}\l( V_{(k/2\ell)\gamma+\Z \gamma},
     \eL (\ell,j)\r),
\end{equation}
where $0\leq j\leq \ell$ and $0\leq k<2\ell$.
Then $W_\ell(j,k)$ denotes the space of multiplicity of $V_{(k/2\ell)\gamma+\Z \gamma}$
in $\eL (\ell,j)$.
Thus, viewing $\eL (\ell,j)$ as a $V_{\Z\gamma}\tensor W_\ell$-module,
we have the following decomposition:
$$
  \eL (\ell,j)= \bigoplus_{k=0}^{2\ell -1} V_{(k/2\ell)\gamma +
  \Z\gamma}\tensor W_\ell(j,k) .
$$
It is shown in \cite{DL} that $W_\ell(j,k)=0$ if $j+k \equiv 1 \mod 2$, and so
\begin{equation}\label{eq:LDEC}
  \eL (\ell,j)=
  \begin{cases}
  \displaystyle
    \bigoplus_{k=0}^{\ell -1} V_{(k/\ell)\gamma + \Z\gamma} \tensor W_\ell(j,2k),
    & \text{ if }j \text{ is  even},
    \\
   \displaystyle \bigoplus_{k=0}^{\ell -1}
     V_{((2k+1)/2\ell)\gamma + \Z\gamma}\tensor W_\ell(j,2k+1),
     & \text{ if }j \text{ is  odd}.
  \end{cases}
\end{equation}
The following basic fact is due to \cite{DL} (see also \cite{Li1}).

\begin{prop}\label{prop:6.1}
  (\cite{DL,Li1})\\
  (1) All $W_\ell(j,k)$, $0\leq j\leq \ell$, $0\leq k \leq
  2\ell-1$, $j\equiv k \mod 2$, are irreducible
  $W_\ell$-modules.\\
  (2) As $W_\ell$-modules, $W_\ell(j_1,k_1)\simeq W_\ell(j_2,k_2)$ if $j_1+j_2=\ell$ and
  $k_2\equiv k_1+\ell \mod 2\ell$.
\end{prop}

We will use $W_\ell$-modules $W_\ell(j,k)$ to study $M_{E_6}$ and $M_{E_7}$.
In \cite{LYY2}, they recursively computed the vacuum characters of these modules.
In the below we shall use the vacuum characters of $W_\ell$-modules
obtained in (loc.\ cit.) without any comments.

\subsection{$V_{\sqrt{2}A_N}$ and $W_{N+1}$-algebra}

Let us consider a lattice VOA $V_{\sqrt{2}A_N}$.
Let $s_{A_N}$ and $\tilde{\w}_{A_N}$ be conformal vectors of $V_{\sqrt{2}A_N}$
defined as in \eqref{eq:2.2} and \eqref{eq:2.3}, respectively.
Then one has an orthogonal decomposition $\w=s_{A_N}+\tilde{\w}_{A_N}$ of
the Virasoro vector $\w$ of $V_{\sqrt{2}A_N}$.
We use the following result established in \cite{LY2}.

\begin{lem}\label{lem:6.2}
  (\cite{LY2})
  There is a VOA-isomorphism $\com_{V_{\sqrt{2}A_N}}(M_{A_N})\simeq W_{N+1}$.
\end{lem}

Due to the lemma above, we shall identify the commutant subalgebra
$\com_{V_{\sqrt{2}A_N}}(M_{A_N})$ with $W_{N+1}$.
By the orthogonality, $(W_{N+1},\tilde{\w}_{A_N})$ forms a sub VOA of $V_{\sqrt{2}A_N}$.
Therefore, $V_{\sqrt{2}A_N}$ contains a full subalgebra isomorphic to
$M_{A_N}\tensor W_{N+1}$.
It is shown in \cite{DLMN} that $M_{A_N}$ contains a full sub VOA isomorphic to
$L(c_1,0)\tensor \cds \tensor L(c_{N},0)$, where $c_m$ denotes the central charge
of the unitary series of the Virasoro algebra:
\begin{equation}\label{eq:6.4}
  c_m:=1-\fr{6}{(m+2)(m+3)},\q m=1,2,\dots.
\end{equation}
By \cite{W}, the irreducible modules over $L(c_m,0)$ are given by the irreducible
highest weight modules $L(c_m,h_{r,s}^m)$ whose highest weights are parameterized as follows.
\begin{equation}\label{eq:6.5}
  h_{r,s}^m:=\dfr{\{ r(m+3)-s(m+2)\}^2-1}{4(m+2)(m+3)},\q
  1\leq r \leq m+1,\q
  1\leq s \leq m+2.
\end{equation}
In \cite{LY2}, the following decomposition of $V_{\sqrt{2}A_N}$ as an $L(c_1,0)\tensor \cds
\tensor L(c_N,0) \tensor W_{N+1}$-module is obtained.
\begin{equation}\label{eq:6.6}
  V_{\sqrt{2} A_N}
  = \bigoplus_{0\leq k_j\leq j+1 \atop {j=0,\dots,N \atop k_j\equiv 0 \mod 2}}
    L(c_1,h^1_{k_0+1,k_1+1})\tensor \cds \tensor L(c_N,h^N_{k_{N-1}+1,k_N+1})
    \tensor W_{N+1}(k_N,0).
\end{equation}
To describe $V_{\sqrt{2} A_N}$ as an $M_{A_N}\tensor W_{N+1}$-module, we introduce
the following notation.
\begin{equation}\label{eq:6.7}
  M_{A_N}(2s):=\bigoplus_{0\leq k_j \leq j+1 \atop {j=0,\dots,N-1\atop k_j\equiv 0 \mod 2}}
  L(c_1,h^1_{k_0+1,k_1+1})\tensor \cds \tensor L(c_N,h^N_{k_{N-1}+1,2s+1}).
\end{equation}
By \eqref{eq:6.6}, we have
\begin{equation}\label{eq:6.8}
  V_{\sqrt{2}A_N}=\bigoplus_{0\leq 2s\leq N+1} M_{A_N}(2s)\tensor W_{N+1}(2s,0)
\end{equation}
and it is shown in \cite{LS} that $M_{A_N}(2s)$, $0\leq 2s\leq N+1$, are inequivalent
irreducible $M_{A_N}$-modules.
By construction, the vacuum characters of $M_{A_N}(2s)$ are obvious.

\subsection{Decomposition of $M_{E_7}$}

Let us recall that the root lattice of type $E_7$ can be written as
$$
  E_7=\left\{(x_1,\ldots,x_8) \in \Q^8 \, \left| \, {\,\mbox{all $x_i$ are
  in}~\mathbb{Z}~ \mbox{or all $x_i$ are
  in}~\frac{1}{2}+\mathbb{Z},}\atop
  {x_1+x_2+\cdots+x_7+x_8=0}\right.\right\}.
$$
Let $\epsilon_i$ be the vector of $\Q^8$ such that the $i$-th entry is $1$
and all the other entries are zero, and set
$$
  N := \Span_\Z\{ -\epsilon_1+\epsilon_2, \dots,-\epsilon_7+\epsilon_8\} .
$$
Then $N$ is a root lattice of type $A_7$.
Let
$\xi=(\shf,\shf,\shf,\shf,-\shf,-\shf,-\shf,-\shf)\in \Q^8$.
Then the root lattice of type $E_7$ can be written as $ E_7=N\cup (\xi+N). $

Let $s_{A_7}$ and $\tilde{\om}_{A_7}=\om -s_{A_7}$ be conformal
vectors defined as in \eqref{eq:2.2} and \eqref{eq:2.3},
respectively. Then $W_8$ is isomorphic to the parafermion algebra
of central charge $7/5$ (cf.\ \cite{LY2, ZF}).
As we have seen, $V_{\sqrt{2}N}$ contains a full sub VOA isomorphic to
$M_{A_7} \tensor W_8$.
The following decomposition is obtained in \cite{LY2,LS}.

\begin{lem}\label{lem:6.3}
  (\cite{LY2,LS})
  As a module over $M_{A_7}\tensor W_8$,
  $$
    V_{\sqrt{2}\xi+\sqrt{2}N}
    \simeq \bigoplus_{0 \leq 2s \leq 8} M_{A_7}(2s)\otimes W_{8}(2s,8).
  $$
  Hence, we have
  $$
    V_{\sqrt{2}E_7}
    = V_{\sqrt{2}N}\oplus V_{\sqrt{2}\xi+\sqrt{2}N}
   \simeq \bigoplus_{0 \leq 2s \leq 8} M_{A_7}(2s)\otimes \left( W_8(2s,0)\oplus W_8(2s,8)\right) .
  $$
\end{lem}

Let $\tilde{\om}'_{E_7}=\tilde{\om}_{A_7}- \tilde{\om}_{E_7}$.
Then $\tilde{\om}_{E_7}$ and $\tilde{\om}'_{E_7}$ are mutually
orthogonal conformal vectors with central charge $7/10$.
Denote
\[
  U:=\com_{V_{\sqrt{2}E_7}}(\vir(s_{A_7}))   \simeq W_{8}(0,0)\oplus W_{8}(0,8).
\]
Then $U$ contains $\tilde{\w}_{E_7}$ and $\tilde{\w}'_{E_7}$ so that
$\vir(\tilde{\w}'_{E_7})\tensor \vir(\tilde{\w}_{E_7})$ is a full sub VOA
of $U$ isomorphic to $L(\sfr{7}{10},0)\otimes L(\sfr{7}{10},0)$.
By the vacuum characters, it is easy to show that
\[
  U\simeq L(\fr{7}{10},0)\otimes L(\fr{7}{10},0)\oplus
  L(\fr{7}{10},\fr{3}{2})\otimes L(\fr{7}{10},\fr{3}{2})
\]
as a $\vir(\tilde{\w}'_{E_7})\tensor \vir(\tilde{\w}_{E_7})$-module.
Moreover, we have

\begin{lem}\label{lem:6.4}
  The lattice VOA $V_{\sqrt{2}E_7}$ can be decomposed as follows:
  \[
     V_{\sqrt{2}E_7} \simeq \bigoplus_{0 \leq 2s \leq 8} M_{A_7}(2s)\otimes U(2s),
  \]
  where
  \[
  \begin{split}
    U(0)\simeq U(8)
    & \simeq L(\frac{7}{10},0)\otimes L(\frac{7}{10},0)
      \oplus L(\frac{7}{10},\frac{3}{2})\otimes L(\frac{7}{10},\frac{3}{2}),
      \\
    U(2)\simeq U(6)
    & \simeq L(\frac{7}{10},\frac{3}{5})\otimes L(\frac{7}{10},\frac{3}{5})
      \oplus L(\frac{7}{10},\frac{1}{10})\otimes L(\frac{7}{10},\frac{1}{10}),
      \\
    U(4)
    & \simeq L(\frac{7}{10},0)\otimes L(\frac{7}{10},\frac{3}{5})\oplus
      L(\frac{7}{10},\frac{3}{2})\otimes L(\frac{7}{10},\frac{1}{10})
      \\
    & \ \oplus L(\frac{7}{10},\frac{3}{5})\otimes L(\frac{7}{10},0)
      \oplus L(\frac{7}{10},\frac{1}{10})\otimes L(\frac{7}{10},\frac{3}{2}).
  \end{split}
  \]
\end{lem}

Since $\w=s_{A_7}+\tilde{\w}'_{E_7}+\tilde{\w}_{E_7}$ is an orthogonal sum,
we have the following isomorphism as $M_{A_7}\tensor \vir(\tilde{\w}'_{E_7})$-modules:
$$
  M_{E_7}
  \simeq M_{A_7}(0) \otimes L(\frac{7}{10}, 0)
  \oplus M_{A_7}(4) \otimes   L(\frac{7}{10}, \frac{3}{5})
  \oplus M_{A_7}(8) \otimes L(\frac{7}{10},0).
$$

\begin{prop}\label{prop:6.5}
  $M_{E_7}$ is generated by its weight $2$ subspace.
\end{prop}

\pf
It is shown in Proposition 5.6 of \cite{LS} that $M_{A_N}$ is generated
by its weight two subspace as a VOA.
Thus so is $M_{A_7}(0)\otimes L(\sfr{7}{10},0)$.
Recall that
\[
  M_{A_7}(2s)\simeq
  \bigoplus_{ { 0\leq k_j\leq j+1,}
  \atop { {j=0,\dots, 7}\atop{\ k_j \equiv0\,\mathrm{mod}\, 2}}}
  L(c_1,h^1_{k_{0}+1,k_1+1})\otimes \cdots \otimes L(c_7,h^7_{k_{6}+1,2s+1}).
\]
Take $(2k_0+1, 2k_1+1, \dots, 2k_6+1,2s+1)=(1,3,3,5,5,7,7,9)$.
Then we obtains a highest weight vector of weight
\[
  (\frac{1}{2},\frac{1}{10},\frac{2}{5},
  \frac{1}{7},\frac{5}{14},\frac{1}{6},\frac{1}{3})
\]
in $M_{A_7}(8)$. Similarly, by taking
$(2k_0+1, 2k_1+1, \dots, 2k_6+1,2s+1)=(1,1,1,1,1,3,5,5)$, $M_{A_7}(4)$ contains
a highest weight vector of weight
\[
  (0,0,0,0,\frac{3}{4},\frac{7}{12},\frac{1}{15}).
\]
Therefore, both $M_{A_7}(4)\otimes L(\sfr{7}{10},\sfr{3}{5})$ and
$M_{A_7}(8)\otimes L(\sfr{7}{10},0)$ contain weight $2$ elements.
Since $M_{A_7}(4)\otimes L(\sfr{7}{10},\sfr{3}{5})$ and
$M_{A_7}(8)\otimes L(\sfr{7}{10},0)$ are irreducible
$M_{A_7}\otimes L(\sfr{7}{10},0)$-modules (cf. \cite{LS}),
$M_{E_7}$ is generated by its weight $2$ subspace.
\qed

\subsection{Decomposition of $M_{E_6}$}

Let us recall that
$$
  E_6=\left\{(x_1,\ldots,x_8)\in \Q^8\, \left| \,
  {\,\mbox{all $x_i$ are in}~\mathbb{Z}~
  \mbox{or all $x_i$ are in}~\frac{1}{2}+\mathbb{Z},}\atop {\mbox{and}~
  x_1+x_8=x_2+\cdots+x_7=0}\right.\right\} .
$$
Define
\[
\begin{split}
  L_1
  &= \l\{(0;x_2, \dots,x_7;0) \in \Z^8 \mid x_2+\ldots+x_7=0 \r\},
  \\
  L_2
  &= \l\{(x_1;0,\dots, 0;x_8) \in \Z^8 \mid x_1+x_8=0 \r\}.
\end{split}
\]
Then $L_1\simeq A_5$, $L_2\simeq A_1$ and it gives an embedding of
$A_5\oplus A_1$ into $E_6$. Set
\[
  \xi=(\frac{1}{2};\frac{1}{2},\frac{1}{2},\frac{1}{2},-\frac{1}{2},
  -\frac{1}{2},-\frac{1}{2};-\frac{1}{2})\quad  \text{ and } \quad
  L=L_1\oplus L_2.
\]
Then we have $E_6=L\cup (\xi+L)$ and
$V_{\sqrt{2}E_6}=V_{\sqrt{2}L}\oplus V_{\sqrt{2}\xi+\sqrt{2}L}$.
Note that
\[
  V_{\sqrt{2}L}
  \simeq V_{\sqrt{2}A_5}\otimes V_{\sqrt{2}A_1}
  \quad
  \text{ and }
  \quad
  V_{\sqrt{2}\xi+\sqrt{2}L}
  \simeq V_{\sqrt{2}\xi_1+\sqrt{2}A_5}\otimes V_{\sqrt{2}\xi_2+\sqrt{2}A_1},
\]
where we have set
$\xi_1=(0;\shf,\shf,\shf,-\shf,-\shf,-\shf;0)$ and
$\xi_2=(\shf;0,0,0,0,0,0;-\shf)$.

Define $s_{A_5}$, $\tilde{\w}_{A_5}\in V_{L_1}$ and
$w^\pm(2\xi_5)\in V_{L_2}$ as in \eqref{eq:2.2}, \eqref{eq:2.3}
and \eqref{eq:4.3}, respectively, and we identify
$(M_{A_5},s_{A_5})$ as a sub VOA of $V_{L_1}\simeq
V_{\sqrt{2}A_5}$. Set
\begin{equation}\label{eq:6.9}
  \w^1:=\tilde{\w}_{A_5}+w^+(2\xi_2)-\tilde{\w}_{E_6}\in V_{\sqrt{2}E_6}\q \mbox{and}\q
  \w^2:= w^-(2\xi_2) \in V_{L_2}.
\end{equation}
Then $\w^1$ and $\w^2$ are conformal vectors of $V_{\sqrt{2}E_6}$
with central charges 25/28 and 1/2, respectively. We also note
that $\vir(\w^1)\simeq L(\sfr{25}{28},0)$ and $\vir(\w^2)\simeq
L(\sfr{1}{2},0)$. One can directly check that
$\w=s_{A_5}+\w^1+\w^2+\tilde{\w}_{E_6}$ is an orthogonal sum (cf.\
Lemma \ref{lem:2.1}). Therefore, $M_{\sqrt{2}E_6}$ contains a full
sub VOA $M_{A_5}\tensor \vir (\w^1)\tensor \vir (\w^2)\simeq
M_{A_5}\tensor L(\sfr{25}{28},0) \tensor L(\sfr{1}{2},0)$. By a
similar method as in \cite{LS, LY2}, one can establish the
following.

\begin{lem}\label{lem:6.6}
  As a module over $M_{A_{5}}\tensor \vir(\w^1)\tensor \vir(\w^2)$,
  $$
  \begin{array}{lll}
    M_{E_6}
    & \simeq
    & M_{A_5}(0)\tensor \l\{ L(\dfr{25}{28},0)\tensor L(\dfr{1}{2},0)
      \oplus L(\dfr{25}{28},\dfr{15}{2})\tensor L(\dfr{1}{2},\dfr{1}{2})\r\}
    \vsb\\
    && \oplus M_{A_5}(2)\tensor \l\{ L(\dfr{25}{28},\dfr{13}{4})\tensor L(\dfr{1}{2},0)
       \oplus L(\dfr{25}{28},\dfr{3}{4})\tensor L(\dfr{1}{2},\dfr{1}{2})\r\}
    \vsb\\
    && \oplus M_{A_5}(4)\tensor \l\{ L(\dfr{25}{28},\dfr{3}{4})\tensor L(\dfr{1}{2},0)
       \oplus L(\dfr{25}{28},\dfr{13}{4})\tensor L(\dfr{1}{2},\dfr{1}{2})\r\}
    \vsb\\
    && \oplus M_{A_5}(6)\tensor \l\{ L(\dfr{25}{28},\dfr{15}{2})\tensor L(\dfr{1}{2},0)
       \oplus L(\dfr{25}{28},0)\tensor L(\dfr{1}{2},\dfr{1}{2})\r\} .
  \end{array}
  $$
\end{lem}

By \eqref{eq:6.7}, $M_{A_5}$ contains a sub VOA isomorphic to
$M_{A_4}$ such that $\com_{M_{A_5}}(M_{A_4})\simeq
L(\sfr{25}{28},0)$. Therefore, $M_{A_5}$ has a full sub VOA
isomorphic to $M_{A_4}\tensor L(\sfr{25}{28},0)$ and we have the
following decomposition:
\begin{equation}\label{eq:6.10}
  M_{A_5}(s)=\bigoplus_{m=0,2,4} M_{A_4}(m)\tensor L(\dfr{25}{28},h^5_{m+1,s+1}).
\end{equation}
Denote the conformal vectors of $M_{A_4}$ and $\com_{M_{A_5}}(M_{A_4})$ by
$u$ and $v$, respectively.
Then $\tilde{\w}_{A_5}=u+v$ is an orthogonal sum and we have the following sequence
of the full sub VOAs of $V_{\sqrt{2}E_6}$:
$$
  M_{A_4}\tensor \vir(v)\tensor \vir(\w^1)\tensor \vir(\w^2)
  \subset
  M_{A_5}\tensor \vir(\w^1)\tensor \tensor \vir(\w^2)
  \subset
  M_{E_6}.
$$
Thus $\com_{M_{E_6}}(M_{A_4})$ contains a full sub VOA $\vir(v)\tensor \vir(\w^1)
\tensor \vir(\w^2)$ isomorphic to $L(\sfr{25}{28},0)\tensor L(\sfr{25}{28},0)
\tensor L(\sfr{1}{2},0)$.
By Lemma \ref{lem:6.6} and \eqref{eq:6.10}, we have

\begin{lem}\label{lem:6.7}
  As a module over $\vir(v)\tensor \vir(\w^1)\tensor \vir(\w^2)$,
  $$
  \begin{array}{l}
    \com_{M_{E_6}}(M_{A_4})
    \vsb\\
    \simeq
     L(\dfr{25}{28},0)\tensor L(\dfr{25}{28},0)\tensor L(\dfr{1}{2},0)
      \oplus L(\dfr{25}{28},\dfr{3}{4})\tensor L(\dfr{25}{28},\dfr{13}{4})
      \tensor L(\dfr{1}{2},0)
    \vsb\\
    \q \oplus L(\dfr{25}{28},\dfr{13}{4})\tensor L(\dfr{25}{28},\dfr{3}{4})
     \tensor L(\dfr{1}{2},0)
     \oplus L(\dfr{25}{28},\dfr{15}{2})\tensor L(\dfr{25}{28},\dfr{15}{2})
     \tensor L(\dfr{1}{2},0)
   \vsb\\
   \q \oplus L(\dfr{25}{28},0)\tensor L(\dfr{25}{28},\dfr{15}{2})
      \tensor L(\dfr{1}{2},\dfr{1}{2})
      \oplus L(\dfr{25}{28},\dfr{3}{4})\tensor L(\dfr{25}{28},\dfr{3}{4})
      \tensor L(\dfr{1}{2},\dfr{1}{2})
   \vsb\\
   \q \oplus L(\dfr{25}{28},\dfr{13}{4})\tensor L(\dfr{25}{28},\dfr{13}{4})
      \tensor L(\dfr{1}{2},\dfr{1}{2})
      \oplus L(\dfr{25}{28},\dfr{15}{2})\tensor L(\dfr{25}{28},0)
      \tensor L(\dfr{1}{2},\dfr{1}{2}).
  \end{array}
  $$
\end{lem}

The key observation is that $\com_{M_{E_6}}(M_{A_4})$ is isomorphic to
a $\Z_2$-orbifold subalgebra of the 5A-algebra considered in \cite{LYY2}.
We refer the proof of the following fact to \cite{LYY2}.

\begin{lem}\label{lem:6.8}(\cite{LYY2})
  $\com_{M_{E_6}}(M_{A_4})$ is generated by its weight two subspace.
\end{lem}

Now we can show the following.

\begin{prop}\label{prop:6.9}
  $M_{E_6}$ is generated by its weight two subspace.
\end{prop}

\pf
Let $U$ be the sub VOA of $M_{E_6}$ generated by the weight two subspace of $M_{E_6}$.
Since $M_{A_5}$ is generated by its weight two subspace (cf.\ Proposition 5.6 of
\cite{LS}), $U$ contains $M_{A_5}$.
It is also shown in \cite{LS} that all $M_{A_5}(s)$, $s=0,2,4,6$, are irreducible
$M_{A_5}$-submodules of $M_{E_6}$.
On the other hand, thanks to Lemma \ref{lem:6.8}, $U$ also contains all the highest
weight vectors for $\vir(\w^1)\tensor \vir(\w^2)$ which appear in the decomposition
in Lemma \ref{lem:6.6}.
Therefore, $U$ contains all the components given in Lemma \ref{lem:6.6} and
hence $U=M_{E_6}$.
\qed

\end{document}